\documentclass[a4paper,11pt]{article}

\usepackage{a4wide}
\usepackage{amsthm,amsmath,amssymb}
\usepackage[all]{xy}
\usepackage{enumerate}
\usepackage{graphics}

\newtheorem{theorem}{Theorem}
\newtheorem{proposition}[theorem]{Proposition}
\newtheorem{lemma}{Lemma}[section]
\theoremstyle{definition}
\newtheorem{remark}[lemma]{Remark}
\newtheorem{definition}[lemma]{Definition}

\DeclareMathOperator{\C}{\mathbb{C}}

\DeclareMathOperator{\Z}{\mathbb{Z}}

\DeclareMathOperator{\OF}{\mathcal{O}}

\def\cdt{\!\cdot\!}
\def\Hom{{\rm Hom}}
\def\Ind{{\rm Ind}}

\begin{document}
\title{Test vectors for trilinear forms when at least one representation is not supercuspidal}

\author{Mladen Dimitrov \thanks{Institut de Math\'ematiques de Jussieu, Universit\'e Paris 7, UFR Math\'ematiques Site Chevaleret, Case 7012, 75205 Paris Cedex 13, France  \texttt{dimitrov@math.jussieu.fr}} 
\and  Louise Nyssen  \thanks{Institut de Math\'ematiques et de Mod\'elisation de Montpellier, Universit\'e Montpellier 2, CC 051, 
Place Eug\`ene Bataillon, 34095 Montpellier Cedex \texttt{lnyssen@math.univ-montp2.fr}} }

\maketitle

\begin{abstract}
Given three irreducible, admissible, infinite dimensional  complex
representations of  $\mathrm{GL}_2(F)$, with $F$ a local  field, the space 
of  trilinear functionals invariant by the group has dimension at most one. 
When it is one we provide an explicit vector on which the 
functional does not vanish assuming that not all three representations are 
supercuspidal.
\end{abstract}


\section{Introduction}

\subsection{What is a test vector? }

Let $F$ be a local non-Archimedean field  with ring of integers $\OF$,
uniformizing parameter $\pi$ and finite residue field. Let  $V_1$, $V_2$ and
$V_3$ be three irreducible, admissible, infinite dimensional  complex
representations of  $G=\mathrm{GL}_2(F)$ with central characters 
$\omega_1$, $\omega_2$ and $\omega_3$  and conductors $n_1$, $n_2$ and $n_3$. 
Using the theory of Gelfand pairs, Dipendra Prasad proves in \cite{P}
that the space of $G$-invariant linear forms on $V_1\otimes V_2
\otimes V_3$, with $G$ acting diagonally,  has dimension at most one and gives 
a precise criterion for this dimension to be  one, that we will now explain.

Let $D^\times$ be the group of invertible elements of the unique quaternion
division algebra  $D$ over $F$, and  denote by  $R$  its unique maximal order.
When $V_i$ is a discrete series representation of $G$, denote
by $V^D_i$  the irreducible representation of $D^\times$ associated  
to $V_i$ by the Jacquet-Langlands correspondence. Again, by
the theory of Gelfand pairs,  
the space of $D^\times$-invariant linear forms on $V^D_1\otimes V^D_2 \otimes
V^D_3$ has dimension at most one. 

A necessary condition for the existence of a  non-zero 
$G$-invariant linear form on $V_1\otimes V_2\otimes V_3$ or a 
non-zero $D^\times$-invariant linear form on $V^D_1\otimes V^D_2 \otimes
V^D_3$, that we will {\it always assume}, is that 
\begin{equation}\label{central}
 \omega_1\omega_2\omega_3=1.
 \end{equation}

\begin{theorem} (\cite[Theorem 1.4]{P},\cite[Theorem 2]{P2})  \label{prasad}
Let $\epsilon(V_1 \otimes V_2 \otimes V_3)=\pm 1$ denote the 
root number of the corresponding $8$-dimensional symplectic representation of the Weil-Deligne group of $F$.  
When all the  $V_i$'s are supercuspidal, 
assume either that $F$ has characteristic zero or that its
residue characteristic is odd. 

Then $\epsilon(V_1 \otimes V_2 \otimes V_3)=1$ if,
and only if, there exists a non-zero $G$-invariant linear form $\ell$ on  
$V_1\otimes V_2 \otimes V_3$, and  $\epsilon(V_1 \otimes V_2 \otimes
V_3)=-1$ if, 
and only if,  all the  $V_i$'s are discrete series representations
of $G$ and there exists a non-zero $D^\times$-invariant linear form
$\ell'$ on
$V^D_1\otimes V^D_2 \otimes V^D_3$.
\end{theorem}

Given a non-zero $G$-invariant linear form $\ell$ on $V_1\otimes V_2
\otimes V_3$, our goal is to find a {\it pure tensor} in  $V_1
\otimes V_2 \otimes V_3$ which is not in the kernel of $\ell$. 
We call such a pure tensor a {\it test vector}.  

\bigskip
Let $v_i$ denote a new vector in $V_i$ (see section \ref{nv}). 
The following results are due to  Dipendra Prasad and Benedict Gross. 
They show that tensor products of new vectors can sometimes be test vectors.

\begin{theorem}\label{vt-000}
\begin{enumerate}  
\item (\cite[Theorem 1.3]{P})  If all the $V_i$'s are unramified principal series, 
then  $v_1 \otimes v_2 \otimes v_3$ is a test vector.

\item (\cite[Proposition 6.3]{GP}) 
Suppose that for $1\leq i\leq 3$, $V_i$ is a twist of the Steinberg
representation by an unramified character $\eta_i$. Then 
 
\begin{itemize}
\item either,  $\eta_1\eta_2\eta_3(\pi)=-1$ and   $v_1 \otimes v_2
  \otimes v_3$  is a test vector.  
\item or,  $\eta_1\eta_2\eta_3(\pi)=1$ and the line in
  $V^D_1\otimes V^D_2 \otimes V^D_3$ 
fixed by $R^\times \times R^\times\times R^\times$ is  not in the kernel of $\ell'$. 
\end{itemize} 
\end{enumerate}  
\end{theorem}

However, as mentioned in \cite[Remark 7.5]{GP}, new vectors do not
always yield  test vectors.  Suppose, for example, that  
$V_1$ and $V_2$ are unramified, whereas $V_3$ is ramified, and 
denote by $K=\mathrm{GL}_2(\OF)$ the standard maximal compact subgroup of
$G$.  
Since $v_1$ and $v_2$ are $K$-invariant and $\ell$ is $G$-equivariant,
$v \mapsto  \ell(v_1\otimes v_2\otimes v)$ defines a $K$-invariant
linear form on $V_3$. In the meantime,  $V_3$ and its contragredient are
 ramified, and  therefore the above linear form has to be zero. In particular
$\ell(v_1\otimes v_2\otimes v_3)=0$. 
To go around this obstruction for new vectors to be test vectors, Gross
and Prasad make a  suggestion, which is  the object of our first
result:

\begin{theorem}\label{vt-00n} If $V_1$ and $V_2$  are unramified and 
$V_3$ has conductor $n\geq 1$, then $\gamma^{n}\cdt v_1 \otimes v_2 \otimes v_3$ 
and $v_1 \otimes \gamma^{n}\cdt v_2 \otimes v_3$ are both test vectors,
where 
$\gamma=\begin{pmatrix}\pi^{-1}  & 0 \\ 0 &  1\end{pmatrix}$.   
\end{theorem}

In general we want to exhibit a test vector as an
{\it explicit privileged}  $G$-orbit  inside the 
$G\times G \times G$-orbit  of $v_1 \otimes v_2 \otimes v_3$, 
where $G$ sits diagonally in  $G\times G \times G$. 
Before stating our main result, let us explain a more general and
systematic  approach in the search for test vectors.

\subsection{The tree for $G$} \label{tree}
The vertices of the tree  are in bijection with 
maximal open compact subgroups of $G$ (or equivalently with lattices in $F^2$, up 
to homothetie) and its edges correspond to Iwahori subgroups of $G$, 
each Iwahori being the intersection of the two  maximal compact 
subgroups sitting at the ends of the edge. 
Every Iwahori being endowed with two canonical $(\OF/\pi)^\times$-valued
characters,  choosing one of those characters  amounts to choosing an 
orientation on the corresponding edge. The standard Iwahori subgroup $I=I_1$ corresponds 
to the edge between $K$ and $\gamma K\gamma^{-1}$, and 
changing the orientation on this edge amounts to replacing
the character $\left(\begin{smallmatrix} a  & b
\\ c &  d\end{smallmatrix}\right)\in I \mapsto (d\mod \pi)$ by 
$\left(\begin{smallmatrix} a  & b
\\ c &  d\end{smallmatrix}\right)\in I \mapsto (a\mod \pi)$. 

More generally, for  $n\geq 1$, the $n$-th standard Iwahori subgroup 
$$I_n=\begin{pmatrix} \OF^\times & \OF \\ \pi^n\OF & \OF^\times\end{pmatrix}$$ 
corresponds to the path between $K$ and 
$\gamma^{n}K\gamma^{-n}$, the set of Iwahori subgroups of depth $n$ 
is in bijection with the set of paths of length $n$ on the tree, 
and choosing an orientation on such a path amounts to 
choosing one of the two $(\OF/\pi^n)^\times$-valued characters 
of the corresponding Iwahori.

The new vector  $v_i$ is by definition a non-zero vector in the unique line
of $V_i$ on which $I_{n_i}$ acts by $\left(\begin{smallmatrix} a  & b
\\ c &  d\end{smallmatrix}\right)\mapsto \omega_i(d)$.
Clearly, for every $n\geq 1$, $G$ acts
transitively on the set of oriented paths of length $n$.  
Hence finding a  $G$-orbit inside the 
$G\times G \times G$-orbit of $v_1 \otimes v_2 \otimes v_3$, 
amounts to finding a  $G$-conjugacy class $I'\times I''\times
I'''$ inside the  $G\times G \times G$-conjugacy class of
$I_{n_1}\times I_{n_2}\times I_{n_3}$.

A most natural way of defining such a  $G$-conjugacy class (almost uniquely) 
is by imposing the smallest of the three compact open subgroups to be the intersection of the two others. 

For instance, the test vector $\gamma^{n}\cdt v_1 \otimes v_2\otimes
v_3$  in  Theorem  \ref{vt-00n} corresponds to  the $G$-conjugacy class
of $\gamma^{n}K\gamma^{-n}\times K\times I_n$. The 
 linear form on $V_3$ given by $v \mapsto  \ell(\gamma^{n}\cdt
 v_1\otimes v_2\otimes v)$  is invariant by  $\gamma^{n}K\gamma^{-n} \cap
 K =I_n$, hence belongs to the new line in the contragredient of
 $V_3$.

Visualized on the tree, the condition on the three compact open subgroups means that 
the longest path should be exactly covered by the  two others, as
shown on each of the following two pictures. 
$$
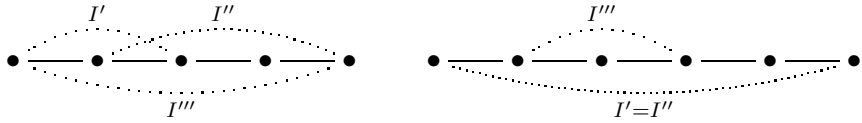
\xymatrix@C=20pt{\bullet \ar@{-}[r]
  \ar@/^1pc/@{.}[rr]^{I'}\ar@/_1pc/@{.}[rrrr]_{I'''} & \bullet
  \ar@{-}[r]\ar@/^1pc/@{.}[rrr]^{I''} & \bullet \ar@{-}[r] & \bullet
  \ar@{-}[r] & \bullet
 &  \bullet \ar@{-}[r]\ar@/_1pc/@{.}[rrrrr]_{I'=I''} & \bullet
 \ar@{-}[r]\ar@/^1pc/@{.}[rr]^{I'''} &\bullet \ar@{-}[r] &\bullet
 \ar@{-}[r] &\bullet \ar@{-}[r] & \bullet}$$

We would like to thank Dipendra Prasad for having shared 
this point of view with us.

\subsection{Main result}

Given an admissible  representation $V$ of $G$ and a character $\eta$
of $F^\times$, we let   $V\otimes\eta$  denote the representation of
$G$ on the same space $V$ with action multiplied  by $\eta\circ\det$, 
called the twist of $V$ by $\eta$.

If $\eta_1$, $\eta_2$ and $\eta_3$ are three characters of $F^\times$
such that $\eta_1 \eta_2 \eta_3 =1$, then the $G$-representations 
$V_1 \otimes V_2 \otimes V_3$ and 
$(V_1 \otimes \eta_1) \otimes(V_2\otimes \eta_2 ) \otimes(V_3\otimes\eta_3)$
are identical, therefore 
\begin{equation}\label{twist}
\Hom_G(V_1 \otimes V_2 \otimes V_3,\C)= 
\Hom_G((V_1 \otimes \eta_1) \otimes(V_2\otimes \eta_2 ) \otimes(V_3\otimes\eta_3),\C).
\end{equation}

Hence finding a test vector in $V_1 \otimes V_2 \otimes V_3$ amounts to 
finding one in $(V_1 \otimes \eta_1) \otimes(V_2\otimes \eta_2 )
\otimes(V_3\otimes\eta_3)$ for some choice of characters  $\eta_1$,
$\eta_2$ and $\eta_3$ such that $\eta_1 \eta_2 \eta_3 =1$.   We would like to
exhibit a test vector in the $G\times G \times G$-orbit of $v'_1
\otimes v'_2 \otimes v'_3$, where $v'_i$ denotes a new vector in
$V_i\otimes \eta_i$, and we want it to be fixed by an  open compact
subgroup as large as possible. Therefore the conductors  of
$V_i\otimes\eta_i$ should be as small as possible. 

Denote by
$n_i^\mathrm{min}$ the minimal possible value for 
the conductor of $V_i\otimes\eta$, when $\eta$ varies. Finally, 
let $n^\mathrm{min}$ denote the minimal possible value  of 
$$\mathrm{cond}(V_1 \otimes \eta_1)+ \mathrm{cond}(V_2\otimes \eta_2
)+\mathrm{cond}(V_3\otimes\eta_3),$$ 
when  $(\eta_1,\eta_2,\eta_3)$ runs over all possible 
triples of characters such that $\eta_1 \eta_2 \eta_3 =1$. 
Note that because of the latter condition, the inequality 
$n^\mathrm{min}\geq n_1^\mathrm{min}+n_2^\mathrm{min}+n_3^\mathrm{min}$ is strict in general. 

Also note that the conductor of  a representation is at least equal to the conductor of it's central character. 
Equality holds if, and only if, the representation is principal and has minimal conductor among it's twists.

\begin{definition}\label{minimal}
\begin{enumerate}
\item The representation  $V_i$ is   {\it minimal} if $n_i=n_i^\mathrm{min}$.
\item The triple of representations $(V_1,V_2,V_3)$ satisfying (\ref{central}) is   {\it minimal} if 
\begin{enumerate}
\item either all non-supercuspidal $V_i$'s are minimal, 
\item or none of the $V_i$'s is supercuspidal and $n^\mathrm{min}=n_1+n_2+n_3$.
\end{enumerate}
\end{enumerate}
\end{definition}

It is clear from the definition that for any 
 $V_1$, $V_2$ and $V_3$, there exist characters  $\eta_1$,
$\eta_2$ and $\eta_3$ such that $\eta_1 \eta_2 \eta_3 =1$
and   $(V_1 \otimes \eta_1,V_2\otimes \eta_2 ,V_3\otimes\eta_3)$ is minimal.
Our main result  states:

\begin{theorem}\label{vt-main}
Suppose that at least one  of $V_1$, $V_2$ and $V_3$ is not
supercuspidal, and that if two amongst them are  supercuspidal with
the same conductor then the third one is a ramified principal series. 
Assume that $(V_1,V_2,V_3)$ is minimal and  $\epsilon(V_1 \otimes V_2 
\otimes V_3)=1$.  If  $n_3\geq n_1$ and $n_3\geq n_2$, then
 $ v_1 \otimes \gamma^{n_3-n_2}\cdot v_2 \otimes v_3$ and  
$\gamma^{n_3-n_1}\cdot  v_1 \otimes v_2 \otimes v_3$ are both test vectors.
 \end{theorem}

\begin{remark}
The test vector $ v_1 \otimes \gamma^{n_3-n_2}\cdot v_2 \otimes v_3$
can be visualized on the tree as follows:  
$$\xymatrix@C=13pt{  K  \ar@{-}[r]
  \ar@/^1pc/@{.}[rrrrr]^{I_{n_1}}\ar@/_1pc/@{.}[rrrrrrr]_{I_{n_3}} & 
\gamma K \gamma^{-1}  \ar@{--}[rr] & &
\gamma^{n_3-n_2} K \gamma^{n_2-n_3}
\ar@{--}[rr]\ar@/^1pc/@{.}[rrrr]^{\gamma^{n_3-n_2}I_{n_2}  
\gamma^{n_2-n_3}} & & \gamma^{n_1} K \gamma^{-n_1} 
  \ar@{--}[rr] & &\gamma^{n_3} K \gamma^{-n_3} }$$ 
\end{remark}

\begin{remark} Assume that $(V_1,V_2,V_3)$ is minimal and that 
at least one  of $V_i$'s is not supercuspidal. Then $\epsilon(V_1 \otimes V_2 \otimes V_3)= -1$ if, and only if, one of the representations, say $V_1$,
is a twist of the Steinberg representation by an unramified character $\eta$ and
$V_2$ is a discrete series whose contragredient is isomorphic  
to  $V_3$ twisted by $\eta$ (see \cite[Propositions 8.4, 8.5, 8.6]{P}).
\end{remark}

\begin{remark} Finding test vectors in the case when  all the $V_i$'s are supercuspidal remains an open question. Consider for example the case
when the $V_i$'s have trivial central characters and share the same conductor $n$.  It is well known that the Atkin-Lehner involution  
$\left(\begin{smallmatrix}0  & 1 \\ \pi^n &  0 \end{smallmatrix}\right)$ acts on 
$v_i$ by the root number $\epsilon(V_i)=\pm 1$. It follows that if $\epsilon(V_1)\epsilon(V_2)\epsilon(V_3)=-1$, then  $\ell(v_1 \otimes v_2\otimes v_3)=0$.

If  $V_1$ is unramified and $V_2$,  $V_3$ are supercuspidal of even conductor $n$, trivial central characters and $\epsilon(V_2)\epsilon(V_3)=-1$, then 
by applying the Atkin-Lehner involution one sees that $\ell(\gamma^{n/2}v_1 \otimes v_2\otimes  v_3)=0$.
Similarly, if  $V_1$ is the Steinberg representation and $V_2$,  $V_3$ are supercuspidal of odd conductor $n$, trivial central characters and $\epsilon(V_2)\epsilon(V_3)=1$, then  
by applying the Atkin-Lehner involution one sees that $\ell(\gamma^{(n-1)/2}v_1 \otimes v_2\otimes  \cdot v_3)=0$.

\end{remark}

\subsection{Application of test vectors to subconvexity}

Test vectors for trilinear forms play an important role in various problems involving  
$L$-functions of triple products of automorphic representations of $\mathrm{GL}(2)$.

One such problem, studied by Bernstein-Reznikov in \cite{BR1,BR2} and
more recently by Michel-Venkatesh in \cite{MV1,MV2}, is about finding 
 {\it subconvexity } bounds for the $L$-functions of automorphic representations of $\mathrm{GL}(2)$
along the critical line. More precisely, given a unitary  automorphic representation $\Pi$ of 
$\mathrm{GL}(N)$ over a number field $E$, the subconvexity bound asserts the existence  of an absolute
constant $\delta>0$ such that  :
$$L(\Pi,1/2)\ll_{E,N}C(\Pi)^{1/4-\delta},$$
where  $C(\Pi)$ denotes the analytic conductor of $\Pi$. 
We refer to \cite{MV2}  for the definition of $C(\Pi)$ and for various applications 
of subconvexity bounds to problems in number theory, such as Hilbert's eleventh problem. 
Let us just mention that the subconvexity bounds follow from the Lindel\"off Conjecture, which is true under the Generalized Riemann Hypothesis.

In \cite[1.2]{MV2} the authors establish the following subconvexity bound for $\mathrm{GL}(2)\times \mathrm{GL}(2)$: 
$$L(\Pi_1\otimes \Pi_2,1/2)\ll_{E,C(\Pi_2)}C(\Pi_1)^{1/2-\delta},$$
and obtain as a corollary subconvexity bounds for $\mathrm{GL}(1)$ and $\mathrm{GL}(2)$.
A key ingredient in their proof is to provide a test vector in the following setup: let 
$F$ be the completion of $E$ at a finite place and denote by $V_i$ the local component of $\Pi_i$ 
at $F$  ($i=1,2$). 
Let $V_3$ be a {\it minimal} principal series representation of $G=\mathrm{GL}_2(F)$ such
that (\ref{central}) is fulfilled, and denote by $\ell$ a normalized $G$-invariant trilinear form on 
$V_1\otimes V_2\otimes V_3$(the process of normalization is explained in \cite[3.4]{MV2}).  Then one needs to find a norm $1$ test vector
$v\otimes v' \otimes v''\in V_1\otimes V_2\otimes V_3$ such that 
$$\ell(v\otimes v' \otimes v'')\gg_{n_2}n_1^{-1/2},$$ 
which can be achieved either by using the test vectors from our main theorem, or 
by a direct computation in the Kirillov model as in \cite[3.6.1]{MV2}.

\subsection{Organization of the paper}
 In section \ref{rappels} we recall
basic facts about  induced admissible representations of $G$ which 
are used in section \ref{Preuve-main} to prove Theorem
\ref{vt-00n} and a slightly more general version of Theorem \ref{vt-main} in the case when at most one of the representations is supercuspidal. Section \ref{2SC} recalls some basic facts about Kirillov models and  contains a proof of Theorem \ref{vt-main} in  the case of two  supercuspidal representations.  Finally, in section \ref{reductible} we
study test vectors in reducible induced representation, 
as initiated in the work of Harris and Scholl \cite{HS}.

\subsection{Acknowledgments}
We would like to  thank Philippe Michel for suggesting the study of
this problem,  
and of course Benedict Gross and Dipendra Prasad for their articles
full of inspiration. 
The second named author would like to  thank also Paul Broussous and
Nicolas Templier for many interesting conversations, and  Wen-Ching
Winnie Li for the opportunity to spend one semester at Penn State
University, where the first draft of this paper was written. Finally,
the first named author would like to thank Dipendra Prasad for
several helpful discussions that took place during his visit of
the Tata Institute of Fundamental Research in december 2008, as well 
as the institution for its hospitality and excellent  conditions for work.

\section{Background on induced admissible representations of
  G}\label{rappels} 

\subsection{New vectors and contragredient representation}\label{nv}

Let $V$ be an irreducible, admissible, infinite dimensional
representation of  $G$ with central character $\omega$.
To  the descending chain of open  compact subgroups of $G$ 
$$ K=I_0 \supset I=I_1 \supset \cdots \supset I_n \supset I_{n+1} \cdots  $$
one can associate an ascending chain of vector spaces for $n\geq 1$:
$$V^{I_{n},\omega}=\left\{v\in V \Big{|}
\begin{pmatrix} a & b \\ c &
      d\end{pmatrix}\cdt v=\omega(d)v \text{ , for all }
\begin{pmatrix} a & b \\ c & d\end{pmatrix}\in I_n \right\}.$$
Put $ V^{I_0,\omega}=V^K$. There exists a minimal $n$ such that the vector space
$V^{I_{n},\omega}$ is non-zero.  
It is necessarily one dimensional, called the new line,  and any
non-zero vector in it is called a {\it new   vector} of $V$ (see \cite{C}).    
The integer $n$ is the {\it conductor} of $V$. 
The representation $V$ is said to be {\it unramified} if $n=0$. 

The contragredient representation $\widetilde{V}$ is the space of smooth linear forms
$\varphi$ on $V$, where $G$  acts as follows: 
$$\forall g\in G, \qquad  \forall v\in V, \qquad (g\cdot \varphi)(v)= \varphi(g^{-1}\cdot v).$$ 

There is an isomorphism $\widetilde{V}\simeq V\otimes \omega^{-1}$, hence
$\widetilde{V}$ and $V$ have the same conductor $n$. Moreover, under
this isomorphism the new line in $\widetilde{V}$ is sent to:
$$\left\{v\in V \Big{|} \begin{pmatrix} a & b \\ c &
      d\end{pmatrix}\cdt v=\omega(a)v \text{ , for all }
\begin{pmatrix} a & b \\ c & d\end{pmatrix}\in I_n \right\},$$
which is the image of the new line in $V$ by the Atkin-Lehner involution 
$\left(\begin{smallmatrix}0  & 1 \\ \pi^n &  0
  \end{smallmatrix}\right)$.

\subsection{Induced representations}\label{notations}
Let $(\rho, W)$ be a smooth representation of a closed subgroup $H$ of
$G$. Let $\Delta_H$ be the modular function on $H$.  
The induction of $\rho$ from $H$ to $G$,  denoted $\Ind_{H}^{G} \rho
$,  is the space   of functions $f$ from $G$ to $W$ 
satisfying the  following two conditions:
\begin{enumerate}
\item  for all  $h \in H$ and  $g \in G$ we have 
$f(hg)=\Delta_H(h)^{-\frac{1}{2}} \rho(h) f(g)$; 
\item  there exists an open compact subgroup $K_f$ of $G$ such that for all 
$k \in K_f$ and  $g \in G$ we have  $f(gk)= f(g)$. 
\end{enumerate}
The action of   $G$ is  by right translation:  for all  
$g, g' \in G, (g\cdot f)(g') = f(g'g)$. 
With the  additional condition  that $f$ must be compactly supported
modulo $H$,  one gets the {\it compact induction} denoted by ${\rm
  ind}_{H}^{G}$.  When $G/H$ is compact, there is no difference between 
$\Ind_{H}^{G}$ and ${\rm ind}_{H}^{G}$.  

\bigskip
Let $B$ be the Borel subgroup of upper triangular matrices in $G$, and
let $T$ be the diagonal torus.  
The character  $\Delta_T$ is trivial and we will use $\Delta_B=\delta^{-1} $ with 
$\delta \begin{pmatrix} a & b \\ 0 & d \end{pmatrix}  = \vert
\frac{a}{d} \vert$  
where $ \vert \enspace \vert$ is the norm on $F$. The
quotient $B \backslash G$ is compact and can be identified with $\mathbb{P}^1(F)$. 

\bigskip

\noindent Let $\mu$ and $\mu'$ be two characters of $F^\times$ and $\chi$ be the character of $B$ given by  
$$ \chi  \begin{pmatrix} a &  * \\ 0 & d \end{pmatrix}  =  \mu(a)\mu'(d).$$
The next two sections are devoted to the study of new vectors in $V
= \Ind_{B}^{G} (\chi)$.

\subsection{New vectors in principal series representations}\label{nv-functions}

Assume that $V = \Ind_{B}^{G} (\chi)$ is  a principal series representation of  $G$,  that is $\mu'\mu^{-1}\neq |\cdot|^{\pm 1}$. 
Then $V$ has  central character $\omega=\mu\mu'$ and conductor $n=\mathrm{cond}(\mu)+\mathrm{cond}(\mu')$. Let $v$ denote a new vector in $V$.

When $V$ is unramified  the function $v:G\rightarrow\C$ is such that for all $b \in B$, $g \in G$ and $k \in K$
$$v(bgk)=\chi(b)\delta(b)^\frac{1}{2} v(g),$$
whereas, if $V$ is ramified, then for all $b \in B$, $g \in G$ and $k= \begin{pmatrix} * &  * \\ * & d  \end{pmatrix} \in I_n$, 
$$v(bgk)=\chi(b)\delta(b)^\frac{1}{2} \omega(d)v(g).$$

We normalize $v$ so that $v(1)=1$ and put 
 $$\alpha^{-1}= \mu(\pi)|\pi|^{\frac{1}{2}} \qquad{\rm and} \qquad \beta^{-1}= \mu'(\pi)|\pi|^{-\frac{1}{2}}.$$

\begin{lemma}\label{calcul-NR} If $V$ is unramified then for all $r\geq 0$, 
$$(\gamma^{r}\cdt v)(k)=
\begin{cases}
\alpha^s\beta^{r-s}  & \text{ , if } k\in I_s\setminus I_{s+1}
\text{  for  } 0\leq s\leq r-1,\\
\alpha^{r} & \text{ , if } k\in I_r.\\
\end{cases}$$

Similarly  for $r\geq 1$, 
\begin{equation*}
(\gamma^{r}\cdot v-\alpha\gamma^{r-1}\cdot v)(k)=
\begin{cases} \alpha^{s}\beta^{r-s}-\alpha^{s+1}\beta^{r-1-s} &
\text{, if } k\in I_s\setminus I_{s+1}, 0\leq s\leq r-1,\\
0 & \text{, if } k\in I_{r},\\
\end{cases}
\end{equation*}
\begin{equation*}
\text{ and  }(\gamma^{r}\cdt v-\beta\gamma^{r-1}\cdt v)(k)=
\begin{cases} \alpha^r(1-\frac{\beta}{\alpha}) &
\text{ , if } k\in I_r,\\
0 & \text{ , if } k\in K \setminus I_{r}.\\
\end{cases} 
\end{equation*}

\end{lemma}

\noindent{\it Proof: }
If $k\in I_r$, then $\gamma^{-r}k\gamma^{r}\in K$, so
$(\gamma^{r}\cdt v)(k)=\alpha^{r}v(\gamma^{-r}k\gamma^{r})= \alpha^{r}$. 
Suppose that $k =\begin{pmatrix}a & b \\ c & d\end{pmatrix}\in
I_s\setminus I_{s+1}$ for some $0\leq s\leq r-1$.
Then $\pi^{-s} c \in {\OF}^{\times}$ and
$$ (\gamma^{r}\cdt v)(k)=\alpha^{r}v\begin{pmatrix}a & \pi^{r}b \\ \pi^{-r}c & d\end{pmatrix} 
= \alpha^{r}v\begin{pmatrix} (ad-bc)\pi^{r-s} & a \\ 0 & \pi^{-r}c \end{pmatrix}
= \alpha^s\beta^{r-s}.$$ 
The second part of the lemma follows by a direct computation. $\hfill \Box $ 

\medskip

For the rest of this section we assume that $V$ is ramified, that is $n \geq 1$. 
We put $$m=\mathrm{cond}(\mu') \qquad {\rm so \quad that} \qquad n-m=\mathrm{cond}(\mu).$$ 

By \cite[pp.305-306]{C} the restriction to $K$ of a new vector 
$v$ is supported by the double coset of 
$\left(\begin{smallmatrix}1 & 0 \\ \pi^{m} & 1  \end{smallmatrix}\right)$ modulo $I_{n}$. 
In particular if  $\mu'$ is unramified ($m=0$), then $v$ is supported by 
$I_{n}\left(\begin{smallmatrix} 1 & 0 \\ 1 & 1  \end{smallmatrix}\right)I_{n}
=I_{n} \left(\begin{smallmatrix} 0 & 1 \\ 1 & 0  \end{smallmatrix}\right)I_{n} =K\setminus I$.

If $1 \leq m \leq n-1$, then $v$ is supported by 
$I_{n} \left(\begin{smallmatrix}1 & 0 \\ \pi^m & 1  \end{smallmatrix}\right)I_{n}
=I_{m} \setminus I_{m+1}$.

If  $\mu$  is unramified, then $v$ is supported by  $I_{n}$. 
We normalize $v$ so that $v \left(\begin{smallmatrix}  1 & 0 \\ \pi^{m} & 1 \end{smallmatrix}\right)=1$.

\begin{lemma}\label{calcul-SP1} 
Suppose that $\mu$ is  unramified and $\mu'$ is ramified. Then, for all $r\geq 0$ and $k\in K$, 
$$ (\gamma^{r}\cdt v)(k)=
\begin{cases}\alpha^r \mu'(d) & \text{ , if } 
k=\begin{pmatrix} * & * \\ * & d  \end{pmatrix}\in I_{{n}+r}, \\
0  & \text{ ,  otherwise}.
\end{cases} $$
$$\Bigl( \gamma^{r}\cdt v-\alpha^{-1} \gamma^{r+1}\cdt v \Bigr)(k)=
\begin{cases}\alpha^r \mu'(d) & \text{ , if } 
k=\begin{pmatrix} * & * \\ * & d  \end{pmatrix}\in I_{n+r}\setminus I_{n+r+1}, \\
0  & \text{ ,  otherwise}.
\end{cases} $$
\end{lemma}

\noindent{\it Proof: } 
For $k=\begin{pmatrix} a & b \\ c & d
\end{pmatrix}\in K$ we have 
$$\alpha^{-r}(\gamma^{r}\cdt v)(k)= v(\gamma^{-r}k\gamma^{r})= v\begin{pmatrix} a &\pi^r b \\ \pi^{-r}c & d  \end{pmatrix}.$$

It is easy to check that for every $s\geq 1$, 
\begin{equation}\label{intersection}
K\cap B\gamma^{r}I_s\gamma^{-r}=I_{s+r}.
\end{equation} 
It follows that $\gamma^{r}\cdt v$ has its support in $I_{n+r}$.
If $k\in I_{n+r}$ then  $c \in {\pi}^{m+r} {\OF}^{\times}$ for some $m \geq n$, $d\in   {\OF}^{\times}$ and  we have the following decomposition:
\begin{equation}\label{decomposition}
\begin{pmatrix} a &\pi^r b \\ \pi^{-r}c & d  \end{pmatrix}=
\begin{pmatrix} \det{k} & \pi^{-m}cb\\ 0 & \pi^{-m-r}cd  \end{pmatrix}
\begin{pmatrix} 1 & 0 \\ \pi^{m} & 1  \end{pmatrix}
\begin{pmatrix} d^{-1} & 0 \\ 0 & \pi^{m+r}c^{-1}  \end{pmatrix}.
\end{equation}

Hence $$\alpha^{-r}(\gamma^{r}\cdt v)(k)=
\mu(\det{k})\mu'(\pi^{-m-r}cd)(\mu\mu')( \pi^{m+r}c^{-1})= \mu'(d).$$ 
$\hfill \Box $

\begin{lemma}\label{calcul-SP2} 
Suppose that $\mu'$ is  unramified and $\mu$ is ramified. 
Then for all $r\geq 0$, 
$$(\gamma^{r}\cdt v)(k)=
\begin{cases}\alpha^s\beta^{r-s}
  \mu\left(\frac{\det{k}}{\pi^{-s}c}\right) & \text{ , if }  
k=\begin{pmatrix} * & * \\ c & *  \end{pmatrix}\in I_s\setminus
I_{s+1}\text{, with } 0\leq s\leq r,  \\
0  & \text{ , if }   k\in I_{r+1}.
\end{cases} $$

Moreover, if  $r \geq 1$, then  
$$\Bigl(\gamma^{r}\cdt v-\beta \gamma^{r-1}\cdt v \Bigr)(k) =
\begin{cases}\alpha^r
  \mu\left(\frac{\det{k}}{\pi^{-r}c}\right) & \text{ , if } k=\begin{pmatrix} * & * \\ c & *  \end{pmatrix}\in I_r\setminus I_{r+1},  \\
0  & \text{ ,  otherwise}.
\end{cases}
$$
\end{lemma}

\noindent{\it Proof: } We follow the pattern of proof of lemma \ref{calcul-SP1}. 
The restriction of $\gamma^{r}\cdt v$ to $K$ is zero outside
$$K\cap B\gamma^{r}(K\setminus I)\gamma^{-r}=K\setminus I_{r+1}.$$ 
For $ 0\leq s\leq r$ and 
$k=\begin{pmatrix} a & b \\ c & d  \end{pmatrix}\in I_s\setminus I_{s+1}$ 
we use the following decomposition:
\begin{equation}\label{decomposition2}
\begin{pmatrix} a &\pi^r b \\ \pi^{-r}c & d  \end{pmatrix}=
\begin{pmatrix} -\frac{\det{k}}{\pi^{-r}c} & a+\frac{\det{k}}{\pi^{-r}c}\\ 0 & \pi^{-r}c 
\end{pmatrix}
\begin{pmatrix} 1 & 0 \\ 1 & 1  \end{pmatrix}
\begin{pmatrix} 1 & 1+\frac{d}{\pi^{-r}c} \\ 0 & -1  \end{pmatrix}.
\end{equation}
Since $d\in\OF$ and $\pi^{r}c^{-1}\in\OF$ we deduce that: 
$$\alpha^{-r}(\gamma^{r}\cdt v)(k)=
\mu\Bigl(\frac{\det{k}}{\pi^{-r}c}\Bigr)
\mu'(-\pi^{-r}c)\left\vert \pi^{r}c^{-1} \right\vert =
\mu\Bigl(\frac{\det{k}}{\pi^{-s}c}\Bigr)\alpha^{s-r}\beta^{r-s}.$$ 
$\hfill \Box $ 

For the sake of completeness, we mention one more result. 
We omit the proof, since it  will not be used in sequel of this paper. 

\begin{lemma}
If $\mu$ and $\mu'$ are both ramified ($0<m<n$), 
then for all $r\geq 0$ and $k\in K$, 
$$(\gamma^{r}\cdt v)(k)=
\begin{cases}\alpha^r \mu\Bigl(\frac{\det{k}}{\pi^{-(m+r)}c}\Bigr)\mu'(d) &
  \text{ , if } 
k=\begin{pmatrix} * & * \\ c & d  \end{pmatrix}\in I_{m+r}\setminus I_{m+r+1}, \\
0  & \text{ ,  otherwise}.
\end{cases}$$
\end{lemma}

\subsection{New vectors in special representations}\label{nv-special}

In this section, we will assume that $\Ind_{B}^{G} (\chi)$ 
is reducible, that is $\frac{\mu'}{\mu}= \vert \cdot \vert^{\pm 1}$.

\subsubsection{Case  $\frac{\mu'}{\mu}= \vert \cdot \vert$}\label{sous-quotient}
In this case, there exists a character $\eta$ of $F^\times$ such that 
$\mu= \eta \vert \cdot \vert^{-\frac{1}{2}}$ and $\mu'= \eta
 \vert \cdot \vert^{\frac{1}{2}}$. 
The representation $\Ind_{B}^{G} ((\eta\circ\det)\delta^{-\frac{1}{2}})$
has length $2$ and has one irreducible one dimensional subspace,   
generated by the function $\eta \circ \det$. 
When $\eta$ is trivial the quotient is called the Steinberg representation, denoted $\mathrm{St}$. 
More generally, the quotient is isomorphic to $\eta \otimes \mathrm{St}$ and is called a special representation. 
There is a short exact sequence 
\begin{equation}\label{speciale-quotient}
0 \rightarrow   \eta \otimes \C \rightarrow  \Ind_{B}^{G} ((\eta\circ\det)
\delta^{-\frac{1}{2}}) \xrightarrow{{\mathrm{proj}}} \eta \otimes \mathrm{St}
\rightarrow 0. 
\end{equation}
The representation $\eta \otimes \mathrm{St}$ is minimal if, 
and only if,  $\eta$ is unramified. 
Then the subspace of $K$-invariant vectors in $\Ind_{B}^{G} ((\eta\circ\det)
\delta^{-\frac{1}{2}})$ is the line $\eta \otimes \C$ with basis
$\eta \circ \det$.  
Since 
$$K = I \, \sqcup \, (B\cap K)\begin{pmatrix}0&1\\1&0\end{pmatrix}I$$ 
there exists  $v^I$ (resp. $v^{K\setminus I}$) in $\Ind_{B}^{G} ((\eta\circ\det)
\delta^{-\frac{1}{2}})$ taking value $1$ (resp. $0$) on $I$ and $0$
(resp. $1$) on $K \setminus I$. 
Both $v^I$ and $v^{K\setminus I}$ are $I$-invariant and  
$v^I + v^{K\setminus I}$ is $K$-invariant. 
Hence  ${\mathrm{proj}}(v^I)=-{\mathrm{proj}}(v^{K\setminus I})$ is
a new vector in  $\eta\otimes \mathrm{St}$ whose conductor is $1$.

Let us compute $\gamma^{r}\cdt v^I$ as a function on $G$. 
As in section \ref{nv-functions}, put 
$$\alpha^{-1}= \mu(\pi)|\pi|^{\frac{1}{2}}= \eta(\pi)  
\qquad{\rm and} \qquad \beta^{-1}= \mu'(\pi)|\pi|^{-\frac{1}{2}}= \eta(\pi).$$
\begin{lemma}\label{calcul-SP-NR1}  For all $r\geq 0$, we have
$(\gamma^{r}\cdt v^I)(k)=
\begin{cases} 
\alpha^r  & \text{ , if } k\in I_{r+1},\\
0 & \text{ , if }k\in K\setminus I_{r+1}.\\
\end{cases}$
\end{lemma}
\noindent{\it Proof: }
By  (\ref{intersection}), we have  $K\cap B \gamma^r I\gamma^{-r}=I_{r+1}$, hence
$\gamma^{r}\cdt v^I$ vanishes on $ K\setminus I_{r+1}$. 

For  $k \in I_{r+1}$,  $\gamma^{-r}k\gamma^r \in I$, hence
$\gamma^{r}\cdt v^I(k)= \alpha^{r} v^I(\gamma^{-r}k\gamma^r)=\alpha^{r}$. 
$\hfill \Box $

\subsubsection{Case  $\frac{\mu'}{\mu}= \vert \cdot \vert^{-1}$}\label{sous-espace}

The notations and results from this section will only be used in 
section \ref{reductible}. 
There exists a character $\eta$ of $F^\times$ such that 
$ \mu= \eta \vert \cdot \vert^{\frac{1}{2}}$ and $ \mu'= \eta\vert \cdot \vert^{-\frac{1}{2}} $.
The representation $\Ind_{B}^{G} ((\eta\circ\det) \delta^{\frac{1}{2}})$
has length $2$ and the special representation $\eta \otimes \mathrm{St}$
is an irreducible subspace of  codimension $1$. There is a short exact
sequence  
\begin{equation}\label{speciale-sous-espace}
0 \rightarrow   \eta \otimes \mathrm{St}  \rightarrow \Ind_{B}^{G} ((\eta\circ\det)
\delta^{\frac{1}{2}})  \xrightarrow{{\mathrm{proj}}^*} \eta \otimes \C \rightarrow 0.
\end{equation}
When $\eta$ is unramified, the space of $K$ invariant vectors in $\Ind_{B}^{G} ((\eta\circ\det) \delta^{\frac{1}{2}})$ is the line generated by the function 
$v^K$ taking  constant value $1$ on $K$, that is  for all $b$ in $B$
and $k$ in $K$: 
$$v^K(bk)=\eta \bigl( \det(b) \bigr) \delta(b) .$$ 
We normalize the linear form ${\mathrm{proj}}^*$ by  ${\mathrm{proj}}^*(v^K)=1$. 
The function $\gamma \cdot v^K - \eta(\pi)^{-1} v^K$, 
whose image by ${\mathrm{proj}}^*$ is $0$, is a new vector in $\eta \otimes \mathrm{St}$.

Let us compute $v^K$ as functions on $G$. 
As in section \ref{nv-functions}, put 
$$\alpha^{-1}= \mu(\pi)|\pi|^{\frac{1}{2}} = \eta(\pi)|\pi| 
\qquad{\rm and} \qquad \beta^{-1}= \mu'(\pi)|\pi|^{-\frac{1}{2}}= \eta(\pi)|\pi|^{-1}$$

\begin{lemma}\label{calcul-SP-NR2}  For all $r\geq 0$, 
$$(\gamma^{r}\cdt v^K)(k)=
\begin{cases} 
\alpha^s\beta^{r-s}  & \text{ , if } k\in I_s\setminus I_{s+1}
\text{  for  } 0\leq s\leq r-1,\\
\alpha^{r} & \text{ , if } k\in I_r.\\
\end{cases}$$

Similarly  for $r\geq 1$, 
\begin{equation*}
(\gamma^{r}\cdt v^K-\alpha\gamma^{r-1}\cdt v^K)(k)=
\begin{cases} \alpha^{s}\beta^{r-s}-\alpha^{s+1}\beta^{r-1-s} & 
\text{, if } k\in I_s\setminus I_{s+1}, 0\leq s\leq r-1,\\
0 & \text{, if } k\in I_{r},\\
\end{cases}
\end{equation*}
\begin{equation*}
\text{ and  }(\gamma^{r}\cdt v^K-\beta\gamma^{r-1}\cdt v^K)(k)=
\begin{cases} \alpha^r(1-\frac{\beta}{\alpha}) & 
\text{ , if } k\in I_r,\\
0 & \text{ , if } k\in K \setminus I_{r}.
\end{cases} 
\end{equation*}
\end{lemma}

It is worth noting that $v^K$ behaves as the  new vector in an
unramified representation (see Lemma \ref{calcul-NR}). The proof is the
same.

\section{The case when at most one representation is supercuspidal}\label{Preuve-main}  

In this section we prove the following result. 
\begin{theorem}\label{vt-01sc}
Assume that $(V_1,V_2,V_3)$ is minimal, $\epsilon(V_1 \otimes V_2 \otimes V_3)=1$ and that at most one representation is supercuspidal. Then, up to a permutation of  the $V_i$'s, exactly one of the following holds: 
\begin{enumerate}[(a)]
\item  $n_3>n_1$, $n_3>n_2$,  and $\gamma^{n_3-n_1}\cdot
  v_1 \otimes v_2 
\otimes v_3$ and $ v_1 \otimes \gamma^{n_3-n_2}\cdot v_2
\otimes v_3$ are both test vectors;
\item $n_1=n_2\geq n_3$,  and  
 $v_1 \otimes v_2\otimes \gamma^i v_3$
is a test vector, for all $0\leq i\leq n_1-n_3$.
\end{enumerate}
\end{theorem}

By  symmetry,  it is enough to  prove in case (a)  that 
$\gamma^{n_3-n_1}\cdt v_1\otimes v_2\otimes v_3$ is a test vector.

\begin{lemma} Under the assumptions in theorem \ref{vt-01sc}, if $n_1$, $n_2$ and $n_3$ are not all equal, then $V_1$ and $V_2$ are non-supercuspidal and minimal. 
\end{lemma}
\noindent{\it Proof: }
Assume first that we are in case (a), that is $n_3>n_1$ and $n_3>n_2$. 
Since all representations of
conductor at most $1$ are  non-supercuspidal and minimal, we may 
assume that  $n_3 \geq 2 $. Moreover by (\ref{central}):
$$\mathrm{cond}(\omega_3)\leq \max(\mathrm{cond}(\omega_1),\mathrm{cond}(\omega_2)) \leq \max(n_1,n_2) < n_3,$$ 
hence $V_3$ is either supercuspidal or non-minimal.
Since $(V_1,V_2,V_3)$ is minimal, this proves our claim in this case. 

Assume next that we are in case (b), that is $n_1=n_2>n_3$. 
As in previous case, we may assume that $n_1=n_2\geq 2$. 
Then if only one amongst $V_1$ and $V_2$ is non-supercuspidal and 
minimal, say $V_1$, one would obtain 
$$\mathrm{cond}(\omega_1)=n_1>\max(n_2-1,n_3)\geq \max(\mathrm{cond}(\omega_2),\mathrm{cond}(\omega_3)),$$
which is false by (\ref{central}). Hence the claim. $\hfill \Box $

If $n_1=n_2=n_3$ then we can assume without loss of generality that
$V_1$ and $V_2$ are non-supercuspidal and minimal. Furthermore, by
Theorem \ref{vt-000} one can assume that the $V_i$'s are not all three
unramified, nor are all three twists of the Steinberg
representation by unramified characters. Finally, if all the three representations have conductor one 
and if exactly one among them is special, we can assume without
loss of generality that this is $V_3$.

\subsection{Choice of models}\label{models}

If $V_i$ is a principal series for $i=1$ or $2$, then by minimality
there exist  characters  $\mu_i$ and $\mu'_i$ of  $F^\times$, 
at least one of which is unramified,  such that
$\mu'_i\mu_i^{-1}\neq |\cdot|^{\pm 1}$ and 
 $$V_i = \Ind_{B}^{G} \chi_i \quad \text{ ,  where} \quad 
 \chi_i  \begin{pmatrix} a & b \\ 0 & d \end{pmatrix}
  =  \mu_i(a)\mu'_i(d).$$
Using the natural isomorphism 
$$\Ind_{B}^{G} \chi_i\cong \Ind_{B}^{G} \chi'_i\text{ , where } \chi'_i  \left(\begin{matrix} a & b \\ 0 & d
  \end{matrix}\right)  =  \mu'_i(a)\mu_i(d)$$  
one can assume that  $\mu_1$ and $\mu'_2$  are  unramified. 

\medskip

If $V_i$ is a special representation, then  by minimality
there exists an  unramified character $\eta_i$ such that
$V_i=\eta_i\otimes \mathrm{St}$. We put then 
$$\mu_i=\eta_i\,\vert\cdot\vert^{-\frac{1}{2}}, \quad
\mu'_i=\eta_i\,\vert\cdot\vert^{\frac{1}{2}} \text{ and }\quad
\chi_i= (\eta_i \circ \det) \delta^{-\frac{1}{2}}  $$
 and choose as  model for $V_i$ the  exact sequence
 (\ref{speciale-quotient}):  
$$0 \rightarrow   \eta_i \otimes \C \rightarrow 
\Ind_{B}^{G} (\chi_i )
\xrightarrow{{\mathrm{proj}_i}}  V_i \rightarrow 0.$$
As new vectors, we choose $v_1 = \mathrm{proj}_1(v_1^I)$ in $V_1$ and 
$v_2 = \mathrm{proj}_2(v_2^{K\setminus I})$ in $V_2$.

\subsection{Going down using Prasad's exact sequence}\label{suites-exactes}

We will now explain how Prasad constructs  a non-zero $G$-invariant linear
form on $V_1\otimes V_2 \otimes V_3$. First, there is a  canonical isomorphism:
\begin{equation}
\Hom_G ( V_1 \otimes V_2 \otimes V_3, \C)\xrightarrow{\sim} 
\Hom_G (V_1 \otimes V_2, \widetilde{V_3}). 
\end{equation}

\begin{lemma}\label{isomorphisme}
We have $$\Hom_G ( V_1 \otimes V_2, \widetilde{V_3})
\xrightarrow{\sim} \Hom_G \Bigl({\rm Res}_{G}\,\Ind_{B \times B}^{G \times 
  G} ( \chi_1 \times \chi_2 ) , \widetilde{V_3}\Bigr),$$
where the restriction is taken with respect to the 
diagonal embedding of $G$ in $G\times G$. 
\end{lemma}

\noindent{\it Proof: } This is clear when $V_1$ and $V_2$ are principal
series. Suppose $V_2=\eta_2\otimes \mathrm{St}$. Tensoring the  exact sequence
(\ref{speciale-quotient}) for $V_2$ with the projective
$G$-module $V_1$   and  taking 
$\Hom_G ( \cdot, \widetilde{V_3} )$ yields a long exact sequence: 
$$ 0 \rightarrow  \Hom_G \Bigl( V_1\otimes V_2, \widetilde{V_3} \Bigr)
         \rightarrow  \Hom_G \Bigl( V_1\otimes \Ind_{B}^{G} (\chi_2) ,
         \widetilde{V_3} \Bigr)   
         \rightarrow  \Hom_G \Bigl( V_1\otimes \eta_2 ,
         \widetilde{V_3} \Bigr). $$
By minimality and by the assumption made in the beginning of 
section \ref{Preuve-main}, we have 
\begin{equation}\label{vanishing}
\Hom_G ( V_1\otimes \eta_2 , \widetilde{V_3})=0.
\end{equation}
Hence there is a canonical isomorphism:
$$\Hom_G \Bigl( V_1\otimes V_2, \widetilde{V_3} \Bigr)
     \xrightarrow{\sim}  \Hom_G \Bigl( V_1\otimes \Ind_{B}^{G}(\chi_2) , \widetilde{V_3} \Bigr).$$
This proves the lemma when $V_1$  is principal
series. Finally, if 
 $V_1=\eta_1\otimes \mathrm{St}$ for some unramified character $\eta_1$, 
then  analogously there  is a canonical isomorphism: 
$$
\Hom_G \Bigl( V_1\otimes \Ind_{B}^{G} (\chi_2) , \widetilde{V_3} \Bigr)
        \xrightarrow{\sim} \Hom_G \Bigl(\Ind_{B}^{G} (\chi_1) \otimes
        \Ind_{B}^{G} (\chi_2) , \widetilde{V_3} \Bigr).$$
$\hfill \Box $

 The action of $G$ on
$(B\times B) \backslash (G\times G) \cong \mathbb{P}^1(F)\times
\mathbb{P}^1(F)$ has precisely two orbits. 
The first is the diagonal $\Delta_{B \backslash G}$, which is closed and
can be identified with $B \backslash G$.  The second   is  its
complement   which is open and  can be identified with $T \backslash
G$ via the bijection: 
 $$ \begin{matrix}
T \backslash G & \longrightarrow & \Bigl( B \backslash G \times B \backslash G \Bigr) \setminus \Delta_{B \backslash G} \\
\hfill Tg & \longmapsto & \left( Bg, B\left(\begin{smallmatrix} 0 & 1
    \\ 1 & 0   \end{smallmatrix}\right)g \right) \hfill 
\end{matrix}$$

Hence, there is  a short exact sequence of $G$-modules:
\begin{equation}\label{courtesuite}
0 \rightarrow {\rm ind}_{T}^{G}( \chi_1\chi'_2 )  \xrightarrow{\mathrm{ext}} 
{\rm Res}_{G}\,\Ind_{B \times B}^{G \times 
  G} ( \chi_1 \times \chi_2 )  \xrightarrow{\mathrm{res}} 
\Ind_{B}^{G} (\chi_1 \chi_2 \delta^{\frac{1}{2}}  )\rightarrow 0.
\end{equation}
The surjection  $\mathrm{res}$ is given by the restriction to the diagonal. 
The injection $\mathrm{ext}$ takes a function 
$h \in {\rm ind}_{T}^{G}( \chi_1\chi'_2 )$ to a function 
$H \in \Ind_{B \times B}^{G \times G} ( \chi_1 \times \chi_2
)$ vanishing on $\Delta_{B \backslash G}$, such that for all $g\in G$ 
$$H \Bigl( g, \begin{pmatrix} 0 & 1 \\ 1 & 0 \end{pmatrix} g
\Bigr) = h(g) \label{rel}.$$

 Applying the functor $\Hom_G \Bigl( \bullet, \widetilde{V_3}
 \Bigr) $ yields  a long exact sequence:  
\begin{multline}\label{longuesuite} 
0 \rightarrow \Hom_G \Bigl( \Ind_{B}^{G} \Bigl(\chi_1 \chi_2
\delta^{\frac{1}{2}}  \Bigr), \widetilde{V_3} \Bigr)  
  \rightarrow  \Hom_G \Bigl( {\rm Res}_{G}\,\Ind_{B \times B}^{G \times 
  G} ( \chi_1 \times \chi_2 ), \widetilde{V_3} \Bigr)  \rightarrow  \\
  \rightarrow  \Hom_G \Bigl( {\rm ind}_{T}^{G}\Bigl( \chi_1\chi'_2
  \Bigr), \widetilde{V_3} \Bigr) \rightarrow 
{\rm Ext}_G^1 \Bigl( \Ind_{B}^{G} \Bigl(\chi_1 \chi_2
\delta^{\frac{1}{2}}  \Bigr), \widetilde{V_3} \Bigr)\rightarrow  
 \cdots  
\end{multline}

\begin{lemma}
 $\Hom_G ( \Ind_{B}^{G} (\chi_1
\chi_2 \delta^{\frac{1}{2}}  ), \widetilde{V_3} ) =0$.
\end{lemma}

\noindent{\it Proof: } If, say  $V_2$ is special, then the claim is exactly
(\ref{vanishing}), so we can assume that $V_1$ and $V_2$ are both principal
series.

Suppose that $\Hom_G ( \Ind_{B}^{G} (\chi_1
\chi_2 \delta^{\frac{1}{2}}  ), \widetilde{V_3} ) \neq 0$, in 
particular,  $V_3$ is not supercuspidal. 

If $V_1$ and $V_2$ are both ramified, this
 contradicts the minimality assumption, namely  that $n^\mathrm{min}=n_1+n_2+n_3$, since 
$n_2=\mathrm{cond}(V_2\otimes \mu_2^{-1})$ whereas $n_3>\mathrm{cond}(V_3\otimes \mu_2)$. 

Otherwise, if for example $V_1$  is unramified, then 
$n_2=n_3>n_1=0$ which is impossible by the assumptions in 
theorem \ref{vt-01sc}.  
$\hfill \Box $

  By \cite[Corollary 5.9]{P} it follows that 
${\rm Ext}_G^1 ( \Ind_{B}^{G}(\chi_1 \chi_2 \delta^{\frac{1}{2}}
), \widetilde{V_3} )=0$, hence (\ref{longuesuite}) yields:
\begin{equation}
\Hom_G \Bigl( {\rm Res}_{G}\,\Ind_{B \times B}^{G \times 
  G} ( \chi_1 \times \chi_2 ), \widetilde{V_3}
\Bigr)\xrightarrow{\sim}
 \Hom_G \Bigl( {\rm ind}_{T}^{G}( \chi_1\chi'_2 ), \widetilde{V_3}
\Bigr).
\end{equation}  
Finally, by  Frobenius reciprocity   
\begin{equation}\Hom_G \Bigl( {\rm ind}_{T}^{G}( \chi_1\chi'_2 ),
  \widetilde{V_3} \Bigr) \xrightarrow{\sim} 
\Hom_T \Bigl( \chi_1\chi'_2  , \widetilde{V_{3\vert T}} \Bigr).
\end{equation}
Since by (\ref{central}) the restriction of $\chi_1\chi'_2$ to the center equals  
$\omega_3^{-1}$, it follows from \cite[Lemmes
8-9]{W} that the latter space is one dimensional. 
Thus, we have five canonically isomorphic lines with corresponding
bases: 
\begin{equation}\label{chaine}
\begin{matrix}
0\neq\ell & \in & {\rm Hom}_G  \Bigl( V_1 \otimes V_2 \otimes V_3 , \C \Bigr) \\
             &     & \downarrow \wr \\   

0\neq\psi & \in & {\rm Hom}_G  \Bigl( \Ind_{B}^{G}(\chi_1) 
\otimes\Ind_{B}^{G}(\chi_2) \otimes V_3 , \C \Bigr) \\
             &     & \downarrow \wr \\   
    
0\neq\Psi     & \in & {\rm Hom}_G \Bigl( {\rm Res}_{G}\,\Ind_{B \times
  B}^{G \times G} ( \chi_1 \times \chi_2 ) ,
\widetilde{V_3} \Bigr) \\     &     & \downarrow \wr \\    
0\neq\Phi    & \in & {\rm Hom}_G \Bigl( {\rm ind}_{T}^{G}(\chi_1\chi'_2), \widetilde{V_3} \Bigr) \\          &     & \downarrow \wr \\    
0\neq\varphi  & \in & {\rm Hom}_T \Bigl( \chi_1\chi'_2  ,
\widetilde{V_{3\vert T}} \Bigr) \\ 
\end{matrix}
\end{equation}

Observe that $\varphi$ is  a linear form on $V_3$ satisfying:
\begin{equation}\label{phi}
\forall t \in T, \qquad \forall v \in V_3, \qquad
\varphi(t\cdt v) = 
(\chi_1\chi'_2)(t)^{-1}\varphi(v).
\end{equation}

Moreover, for all $v\in \Ind_{B}^{G}(\chi_1)$, $v'\in \Ind_{B}^{G}(\chi_2)$
and $v''\in V_3$, we have the formula: 
\begin{equation}
\ell(\mathrm{proj}_1(v)\otimes \mathrm{proj}_2(v')\otimes v''  )=\psi(v\otimes v'\otimes
v'')= \int_{T\backslash G} v(g) v'\Bigl(\left(\begin{smallmatrix} 0 & 1
    \\ 1 & 0   \end{smallmatrix}\right)g\Bigr)\varphi(g\cdot
v'')dg, 
\end{equation}
where for $i=1,2$,  $\mathrm{proj}_i$ is the map defined in (\ref{speciale-quotient}),
if $V_i$ is special, and identity otherwise.

\subsection{Going up}\label{goingup}

\begin{lemma}\label{lemmeV3}
For all $i\in \Z$, $\varphi(\gamma^i\cdt v_3) \neq 0$.
\end{lemma}

\noindent{\it Proof: }
Take any ${v_0}\in V_3$ such
that $\varphi({v_0}) \neq 0$. By smoothness ${v_0}$ is fixed by the 
principal congruence subgroup $\ker(K\rightarrow
\mathrm{GL}_2(\OF/\pi^{s_0}))$, for some ${s_0}\geq 0$. Then 
$\varphi(\gamma^{{s_0}}\cdt {v_0})=(\mu_1\mu'_2)(\pi^{s_0})\varphi({v_0}) \neq
0$ and $\gamma^{{s_0}}\cdt {v_0}$ is fixed by the congruence subgroup 
$$I_{2{s_0}}^{1}:=\left\{k\in K \Big{|} k\equiv \begin{pmatrix} 1 & * \\ 0 &
      1\end{pmatrix} \pmod{\pi^{2{s_0}}} \right\}.$$
Hence $\varphi(V_3^{I_{s}^{1}})\neq \{0\}$, for all $s\geq 2s_0$. 
Since $I_{s}/I_{s}^{1}$ is a finite abelian group, $V_3^{I_{s}^{1}}$
decomposes as a direct sum  of spaces indexed by the characters of  $I_{s}/I_{s}^{1}$. 
By (\ref{phi}) and by the fact that  $\mu_1\mu'_2$ is
unramified, $\varphi$ vanishes on all summands of 
$V_3^{I_{s}^{1}}$ other than  $V_3^{I_{s},\omega_3}$ (defined in
section \ref{nv}). Hence $\varphi(V_3^{I_{s},\omega_3})\neq \{0\}$.
By  \cite[p.306]{C} the space $V_3^{I_{s},\omega_3}$ has the
following basis: 
$$\left(  v_3\quad, \quad \gamma\cdt v_3\quad,  \dots
  ,\quad\gamma^{s-n_3}\cdt v_3 \right).$$ 
It follows that  $\varphi(\gamma^{i}\cdt v_3)\neq 0$ for some $i\in\Z$,
hence  by  (\ref{phi}), $\varphi(\gamma^{i}\cdt v_3)\neq 0$  
for all $i\in\Z$. 

Note that  the claim also follows from the first case in \cite[Proposition
2.6]{GP} applied  to the split torus $T$ of $G$.  $\hfill \Box$

\bigskip

Let  $n=\max(n_1,n_2,n_3)\geq 1$ and put 
$$ {J_n}=\begin{pmatrix} 1 &\OF \\ 
 \pi^{n}\OF& 1\end{pmatrix}. $$ 

Consider the unique function  $h\in {\rm ind}_{T}^{G}( \chi_1\chi'_2 )$
which is zero outside the open compact subset $T {J_n}$ of $T\backslash
G$ and such that for all $b_0\in \OF$
and $ c_0 \in \pi^{n}\OF$ we have 
$h\left(\begin{smallmatrix} 1 & b_0 \\ c_0 &  1\end{smallmatrix}\right) = 1$.

For every $0\leq i\leq n-n_3$,  ${J_n}$ fixes
$\gamma^i\cdt v_3$.

\medskip
By definition, the function   $g  \mapsto  h(g)\varphi( g\cdt v_3)$
 factors through $G\rightarrow T\backslash G$ and
by lemma \ref{lemmeV3}: 
\begin{equation}\label{Phi-h-v_3}
\Bigl( \Phi(h) \Bigr)(\gamma^i\cdt v_3) = 
\int_{T \backslash G} \! h(g) \, \varphi( g\gamma^i\cdt v_3) dg
=\varphi(\gamma^i\cdt v_3 ) \int_{ {J_n}} \! dk_0\neq 0 .
\end{equation}

Now, we will compute $H=\mathrm{ext}(h)$ as a function on $G\times G$.
Recall that $H:G\times G\rightarrow \C$ is the unique  function 
satisfying: 
\begin{enumerate}
\item for all $b_1,b_2\in B$, $g_1,g_2\in G$,
$H(b_1g_1,b_2g_2)=\chi_1(b_1)\chi_2(b_2)
\delta^{\frac{1}{2}}(b_1b_2)H(g_1,g_2)$, 
\item for all $g\in G$,  $H(g,g)=0$ and $H(g,
\left(\begin{smallmatrix}  0 & 1 \\ 1 &
  0 \end{smallmatrix}\right)g) = h(g)$. 
\end{enumerate}
 Since $G=BK$, $H$ is uniquely determined by its restriction to
$K\times K$. 
Following the notations of section \ref{nv-functions} put 
$$\alpha_i^{-1}= \mu_i(\pi)|\pi|^{\frac{1}{2}}\quad \text{and} \qquad 
 \beta_i^{-1}= \mu'_i(\pi)|\pi|^{-\frac{1}{2}}.$$ 

\begin{lemma} \label{FV} 
For  all $k_1 =\begin{pmatrix} * & * \\ c_1 & d_1 \end{pmatrix}$ 
and $k_2 =\begin{pmatrix} * & * \\ c_2 & d_2 \end{pmatrix}$ in $K$ we
have 
\begin{equation}
H(k_1,k_2)= \begin{cases}
\omega_1(d_1)
\omega_2\left(\frac{-\det{k_2}}{c_2}\right)
 & \text{ , if } k_1\in I_n 
\text{ and } k_2\in K\setminus I, 
\\0 & \text{ , otherwise}.
\end{cases} \end{equation}
\end {lemma}

\noindent{\it Proof: } 
By definition $H(k_1,k_2)=0$ unless there exist
$k_0=\begin{pmatrix} 1 & b_0 \\ c_0 & 1  \end{pmatrix}
\in {J_n}$  such that 
$$k_1k_0^{-1}\in B \qquad \text{ and } \qquad
k_2k_0^{-1}\begin{pmatrix} 0 & 1 \\ 1 & 0 
\end{pmatrix}\in B ,$$ in which case
\begin{equation}\label{formula-h1}
H(k_1,k_2)=\chi_1(k_1k_0^{-1})\chi_2 \Bigl(
k_2k_0^{-1}\begin{pmatrix} 0 & 1 \\ 1 & 0 \end{pmatrix} \Bigr)  
\delta^{\frac{1}{2}}\Bigl(k_1k_0^{-1}k_2k_0^{-1}\begin{pmatrix} 0 & 1
  \\ 1 & 0 \end{pmatrix} \Bigr) h(k_0). 
\end{equation}
From $k_1k_0^{-1}\in B$, we deduce that $\frac{c_1}{d_1}=c_0\in \pi^{n}{\OF}$.
From $k_2k_0^{-1}\begin{pmatrix} 0 & 1 \\ 1 & 0 \end{pmatrix}\in B $  
we deduce that $\frac{d_2}{c_2}=b_0\in {\OF}$. 
Since, for $i \in\{1,2\}$,  both $c_i$ and $d_i$ are in $\OF$, 
and at least one is in $\OF^\times$ it follows that 
\begin{equation}\label{formula-h2}
d_1, c_2 \in \OF^\times , \qquad d_2\in \OF \qquad \text{and} \qquad  c_1\in
\pi^{n}{\OF}.
\end{equation}
  
 Hence  $k_1\in  I_{n}$ and   $k_2\in K \setminus I$. Moreover  
$$k_1k_0^{-1}=\begin{pmatrix}\frac{\det{k_1}}{d_1\det{k_0}} & * \\ 0 &
  d_1 \end{pmatrix} \text{ and } 
k_2k_0^{-1}\begin{pmatrix} 0 & 1 \\ 1 & 0\end{pmatrix}
=\begin{pmatrix}  \frac{-\det{k_2}}{c_2\det{k_0}} & * \\ 0 & c_2 \end{pmatrix}.$$
Since $n \geq n_2$ and $n\geq 1$ we have $\mu_2(\det{k_0})=1$, hence
\begin{equation}\label{calculH}
H(k_1,k_2) = \mu'_1(d_1)\mu_2\left(\frac{-\det{k_2}}{c_2}\right) =\omega_1(d_1)
\omega_2\left(\frac{-\det{k_2}}{c_2}\right).
\end{equation}
Conversely, if  $k_1\in  I_{n} $ and $k_2\in K \setminus I$ one can take 
$k_0=\begin{pmatrix} 1  & d_2c_2^{-1} \\ c_1d_1^{-1} & 1  \end{pmatrix}$.
\hfill $\Box$

\begin{remark}
One can define $h$ and compute the corresponding $H$ 
for values of $n$ smaller than $\max(n_1,n_2,n_3)$.  
However, $H$ does not need to  decompose as 
a product of functions of one variable as in the above lemma, and the
 corresponding element in $V_1\otimes V_2$ will not be a pure
tensor. For example, if $n_3=0$ and  $n_1=n_2\geq 0$, we can take $n=0$ 
and put ${J_0}=\left(\begin{smallmatrix} 1 &\OF \\ 
 \OF& 1\end{smallmatrix}\right)\cap \mathrm{GL}_2(F)$. Then by
(\ref{formula-h1}) and (\ref{formula-h2}) one finds that 
 for all $k_1\in K$ and  $k_2\in K$ 
$$H(k_1,k_2)=\begin{cases} \frac{\omega_2(-\det
    k_2)}{\mu_1\mu_2|\cdot|(d_1c_2-c_1d_2)}
& \text{ , if } d_1\in \OF^\times, 
  \enspace c_2\in \OF^\times \text{ and }  d_1c_2\neq c_1d_2; \\
0 & \text{ ,  otherwise}. \end{cases}$$
\end{remark}

\bigskip
Now, we want to express $H \in V_1\otimes V_2$
in terms of the new vectors $v_1$ and $v_2$. 
Put \begin{equation}\label{v12}
\begin{split}
v_1^*=\begin{cases} 
\gamma^{n}\cdt v_1-\beta_1 \gamma^{n-1}\cdt v_1  & \text{ , if }  V_1
\text{ is unramified,}\\
\gamma^{n-1}\cdt v_1^I  & \text{ , if }  V_1
\text{ is special,}\\
\gamma^{n-n_1}\cdt v_1 &  \text{ , otherwise, }
\end{cases}  \\
\text{ and } v_2^*=\begin{cases} 
 v_2-\alpha_2^{-1} \gamma\cdt v_2  & \text{ ,  if } V_2 \text{ is unramified,}\\
v_2^{K \setminus I} & \text{ ,  if } V_2 \text{ is special,}\\
 v_2 &  \text{ , otherwise.}\\ 
\end{cases}
\end{split}
\end{equation}

\begin{lemma} \label{calcul-H}
With the notations of (\ref{v12}), $H$ is a non-zero multiple of $
v_1^* \otimes v_2^*$. 
\end{lemma}
\noindent{\it Proof: } Both $H$ and $ v_1^* \otimes v_2^*$ are elements in 
$\Ind_{B \times B}^{G \times G} \Bigl( \chi_1 \times \chi_2
\Bigr)$, hence it is enough to compare their restrictions to $K\times
K$. By Lemmas \ref{calcul-NR}, \ref{calcul-SP1}, \ref{calcul-SP2},
\ref{calcul-SP-NR1} and  
 \ref{FV} both restrictions are supported by  $I_{n}\times  (K\setminus I)$.  

In order  to avoid repetitions or cumbersome notations, we will only give
the final result: 
\begin{equation}\label{lambda12} 
\begin{split}
 &H=\lambda_1\lambda_2\mu_2(-1) \alpha_1^{n_1-n} (v_1^* \otimes v_2^*) \text{ , where }\\
&\lambda_i=\begin{cases} \Bigl(1-\frac{\beta_i}{\alpha_i}\Bigr)^{-1} &  \text{ , if }  V_i \text{ is unramified,}\\
1  & \text{ , if }  V_i \text{ is ramified. }
\end{cases}  
 \end{split}
\end{equation}
If $V_i$ is unramified ($i=1,2$), then $\beta_i \not=
\alpha_i$ and  $\lambda_i \not=0$.  \hfill $\square$

\bigskip
Since by definition, for any $v \in V_3$, we have 
 $\psi(H\otimes v)=\Psi(H)(v)=\Phi(h)(v)$, it
 follows  from Lemma \ref{calcul-H} and (\ref{Phi-h-v_3}) that
for every $i$, $0\leq i\leq n-n_3$:   
\begin{equation}\label{l-F-v_3} 
\psi(v_1^*\otimes v_2^* \otimes \gamma^i\cdt v_3)\neq 0.
\end{equation}

At this stage, we do have an explicit test vector, which is
$\mathrm{proj}_1(v_1^*) \otimes \mathrm{proj}_2(v_2^*) \otimes v_3\in V_1\otimes
V_2\otimes V_3$. By section \ref{sous-quotient} we have : 
\begin{equation}\label{v12bis}
\begin{split}
\mathrm{proj}_1(v_1^*)=\begin{cases} 
\gamma^{n}\cdt v_1-\beta_1 \gamma^{n-1}\cdt v_1  & \text{ , if }  V_1
\text{ is unramified,}\\
\gamma^{n-n_1}\cdt v_1 &  \text{ , otherwise, }
\end{cases}  \\
\text{ and } \mathrm{proj}_2(v_2^*)=\begin{cases} 
 v_2-\alpha_2^{-1} \gamma\cdt v_2  & \text{ ,  if } V_2 \text{ is unramified,}\\
 v_2 &  \text{ , otherwise.}\\ 
\end{cases}
\end{split}
\end{equation}
In the next two sections we will simplify it and deduce 
Theorems \ref{vt-00n} and \ref{vt-01sc}.

\subsection{Proof of Theorem \ref{vt-00n}}

Suppose that $n_1=n_2=0$, so that  $n=\max(n_1,n_2,n_3)=n_3\geq 1$. Then   
(\ref{l-F-v_3}) yields:  
$$\ell\Bigl((\gamma^{n}\cdt v_1-\beta_1\gamma^{n-1}\cdt v_1)\otimes
(\gamma\cdt v_2-\alpha_2v_2)\otimes v_3\Bigr)\neq 0. $$

This expression can be simplified as follows. Consider for $m\geq 0$ 
the linear form:
$$\psi_m(\bullet)=\ell(\gamma^{m}\cdt v_1\otimes  v_2
\otimes\bullet)\in\widetilde{V_3}.$$  

As observed in the introduction, 
$\psi_m$ is invariant by $\gamma^{m}K\gamma^{-m}\cap K =I_m $, hence vanishes 
for $m<n=\mathrm{cond}(\widetilde{V_3})$. Therefore, for ${n}\geq 2$: 
\begin{equation*}\begin{matrix}
\ell\Bigl((\gamma^{n}\cdt v_1-\beta_1\gamma^{n-1}\cdt v_1)\otimes
(\gamma\cdt v_2-\alpha_2v_2)\otimes v_3\Bigr)\hfill \\
\hskip1cm = -\alpha_2\psi_{n}(v_3)+\beta_1\alpha_2\psi_{{n}-1}(v_3)+\psi_{{n}-1}(\gamma^{-1} \cdt v_3)-\beta_1\psi_{{n}-2}(\gamma^{-1} \cdt v_3)\hfill\\
\hskip1cm =  -\alpha_2\psi_{n}(v_3)\hfill\\
\hskip1cm =  -\alpha_2\ell(\gamma^{n}\cdt v_1\otimes  v_2\otimes v_3) \neq 0.\hfill\\
\end{matrix}\end{equation*}

If ${n} = 1$, only the two terms in the middle vanish and we obtain
$$\alpha_2\ell(\gamma\cdt v_1\otimes  v_2\otimes v_3)
+\beta_1\ell(v_1\otimes \gamma\cdt v_2\otimes v_3)\neq 0.$$

Put 
$w=\begin{pmatrix}0 & 1 \\ \pi & 0\end{pmatrix}$. Then 
$w \gamma = \begin{pmatrix} 0 & 1 \\ 1 & 0  \end{pmatrix}\in K$ 
and $\gamma^{-1} w = \begin{pmatrix} 0 & \pi \\ \pi & 0  \end{pmatrix}\in \pi
K$. Hence: 
$$\begin{matrix}
\beta_1\ell(v_1\otimes \gamma\cdt v_2\otimes v_3)\hfill 
& = & \beta_1\ell(\gamma\gamma^{-1} w\cdt v_1\otimes w\gamma\cdt v_2\otimes
w\cdt v_3)\hfill \\ 
& = & \beta_1\omega_1(\pi)\ell\bigl(\gamma\cdt v_1\otimes  v_2\otimes
w\cdt v_3\bigr)\hfill \\ 
& = & \alpha_1^{-1}\ell(\gamma\cdt v_1\otimes  v_2\otimes w\cdt v_3).\hfill \\
\end{matrix}$$

Therefore 
$$\ell\Bigl(\gamma\cdt v_1\otimes v_2\otimes(w\cdt v_3+\alpha_1\alpha_2v_3)\Bigr)\neq 0.$$ 
In particular 
$$\Psi({\gamma} \cdt v_1\otimes v_2)\neq 0.$$
Since $\gamma K {\gamma} \cap K = I$, 
$\Psi(\gamma\cdt v_1\otimes v_2) \in \widetilde{V_3}^{I,{\omega_3}^{-1}}$,
cannot vanish on the line ${V_3}^{I,{\omega_3}}$, which is generated by $v_3$, and therefore   
$$\ell(\gamma\cdt v_1\otimes v_2\otimes v_3)=\Psi(\gamma\cdt v_1\otimes v_2)(v_3)\neq 0.$$  
Hence, if  $n\geq 1$, $\gamma^{n}\cdt v_1\otimes v_2\otimes v_3$
is a test vector. By symmetry $ v_1\otimes \gamma^{n}\cdt v_2\otimes
v_3$ is a test vector too. This completes the proof of Theorem \ref{vt-00n}.
\hfill $\square$

\subsection{End of the proof of Theorem \ref{vt-01sc}}

By Theorem \ref{vt-00n} we may assume that $V_1$ or  $V_2$ is
ramified.

\medskip
If  $V_1$ and $V_2$ are both ramified then Theorem \ref{vt-01sc}
follows directly from     (\ref{l-F-v_3}) and (\ref{v12bis}). 

\medskip
If $V_1$ is unramified  (\ref{l-F-v_3}) yields: 
$$\ell\Bigl((\gamma^{n}\cdt v_1-\beta_1\gamma^{n-1}\cdt v_1)\otimes
v_2\otimes v_3\Bigr)\neq 0.$$ 

Since  $n_1=0<n_2$, we are in case (a) of Theorem \ref{vt-01sc}, 
hence $n_2<n_3=n$, which implies   $\gamma^{n_3-1}K\gamma^{1-n_3}\cap
I_{n_2} = I_{{n_3}-1}$ and   
$$\ell(\gamma^{n_3-1}\cdt v_1\otimes v_2\otimes
\bullet)\in \widetilde{V_3}^{I_{n_3-1},\omega_3^{-1}}=\{0\}.$$ 

Therefore  $\ell(\gamma^{n_3}\cdt v_1\otimes v_2\otimes v_3)\neq 0$,  
that is $\gamma^{n_3}\cdt v_1\otimes v_2\otimes v_3$ is a test vector. 

\medskip
Finally, if $V_2$  is unramified (\ref{l-F-v_3}) yields:  
$$\ell\Bigl(\gamma^{n_3-n_1}\cdt v_1 \otimes (\gamma\cdt v_2-\alpha_2 v_2)\otimes v_3\Bigr)\neq 0. $$

Since  $n_2=0<n_1$, we are in case (a) of Theorem \ref{vt-01sc}, 
hence $n_1<n_3=n$, which implies 
$$\ell(\gamma^{n_3-n_1-1}\cdt v_1\otimes v_2\otimes \bullet)\in \widetilde{V_3}^{I_{{n_3}-1},\omega_3^{-1}}=\{0\}.$$ 

It follows that $\ell(\gamma^{n_3-n_1}\cdt v_1 \otimes \gamma\cdt v_2\otimes v_3)=
\ell(\gamma^{n_3-n_1-1}\cdt v_1 \otimes v_2\otimes\gamma^{-1} \cdt
v_3)=0$.

Therefore   $\ell(\gamma^{n_3-n_1}\cdt v_1 \otimes v_2\otimes v_3)\neq 0$, 
that is $\gamma^{n_3-n_1}\cdt v_1 \otimes v_2\otimes v_3$ is a test vector. 

The proof of Theorem \ref{vt-01sc} is now complete. \hfill $\square$

\section{Proof of  Theorem \ref{vt-main} when two of the 
representations are supercuspidal}\label{2SC}  

The proof in this case follows the original approach of Prasad \cite[page 18]{P}. We are indebted to Paul Broussous who has first obtained and 
shared with us some of the results described here. 

Suppose given $V_1$, $V_2$ and $V_3$ as in theorem \ref{vt-main}
and such that exactly two of the  $V_i$'s are supercuspidal. 
The condition  (\ref{central}) forces the representation 
with the largest conductor $V_3$ to be supercuspidal
and we may assume that $V_2$ is supercuspidal too, whereas 
$V_1$ is minimal.

\subsection{Kirillov model for supercuspidal representations}\label{kirillov}
Suppose given an irreducible supercuspidal representation $V$ of $G$ 
with central character $\omega$. 
Fix a non-trivial additive character $\psi$ on $F$ of conductor $0$.
We  identify $F$ with the unipotent subgroup $N$ of $B$ and 
denote by $\psi\boxtimes\omega$ the corresponding character of $NF^\times $. 
Then the compactly induced representation 
$\mathrm{ind}_{NF^\times}^B\left(\psi\boxtimes\omega\right)$ is naturally isomorphic to the
space $\mathcal{C}_c^\infty(F^\times)$ of compactly supported
locally constant functions on $F^\times$ on which $B$ acts as follows: 
\begin{equation}\label{kirillov-action}
\begin{pmatrix} a & b \\ 0 & d  \end{pmatrix}\cdot f (x)=
\omega(d)\psi\left(\frac{b}{d}\right)f\left(\frac{ax}{d}\right).
\end{equation}

It is well known (see \cite[\S4.7]{B} that  the restriction of $V$ to
$B$ is irreducible and isomorphic to
$\mathrm{ind}_{NF^\times}^B\left(\psi\boxtimes\omega\right)$. In other
terms  there is an unique way to endow the latter 
with a $G$-action making it isomorphic to $V$. Hence the 
action of  $B$ on  $\mathcal{C}_c^\infty(F^\times)$ defined in 
(\ref{kirillov-action}) can be uniquely extended to $G$ so that the 
resulting representation is isomorphic to $V$. It is called the
{\it Kirillov model} of $V$, with respect to $\psi$. 

The characteristic function of $\mathcal{O}^\times$ is a new vector in the Kirillov model. 

\subsection{Choice of models}
We first choose a model for  $V_1$. Consider the character $\chi_1$ of $B$
defined by $\chi_1\left(\begin{smallmatrix} a & * \\ 0 & d \end{smallmatrix}\right)  = \vert
\frac{a}{d} \vert^{-\frac{1}{2}}  \omega_1(d)$. 
The claim of the theorem is invariant by  unramified twists.  
 By the minimality assumption, after  twisting 
$V_1$ by an appropriate  unramified character (and  $V_2$  by its inverse), 
we can assume  either that $V_1=\Ind_{B}^{G} \chi_1 $, or  that $V_1$
is the Steinberg  representation. In both cases $V_1$ is the unique
irreducible   quotient of $\Ind_{B}^{G} \chi_1 $. 

\begin{lemma} \label{inclusion}
The natural inclusion of  $\widetilde{V_1}$ in $\Ind_{B}^{G} (\chi_1^{-1})$
induces an  isomorphism:
$$\Hom_G ( V_2 \otimes V_3, \widetilde{V_1})
\xrightarrow{\sim} \Hom_G ( V_2 \otimes V_3, \Ind_B^G(\chi_1^{-1})).$$
\end{lemma}
\noindent{\it Proof: }  The lemma is clear if $V_1$ is a principal series.
If $V_1$ is the Steinberg  representation, the condition  
$\epsilon(V_1\otimes V_2 \otimes V_3)=1$ implies that 
$\Hom_G ( V_2 \otimes V_3, \C)=\Hom_G ( V_2, \widetilde{V_3})=0$. 
The lemma then follows from the long exact sequence obtained by
applying the functor $\Hom_G ( V_2 \otimes V_3,\bullet)$ to the 
short exact sequence $(\ref{speciale-sous-espace})$. \hfill $\square$

By Frobenius reciprocity: 
$$\Hom_G ( V_2 \otimes V_3, \Ind_B^G(\chi_1^{-1}))
\xrightarrow{\sim} \Hom_B ( V_2 \otimes V_3, \chi_1^{-1}\delta^{\frac{1}{2}}).$$

Let us  choose Kirillov models for $V_2$ (resp.  $V_3$) with respect
to $\psi$ (resp. $\overline{\psi}$),  
so that vectors in  $V_2$ and  $V_3$ are elements in
$\mathcal{C}_c^\infty(F^\times)$. 
For $v'\in V_2$ and  $v''\in V_3$ we  define:
\begin{equation}
\Phi(v',v'')=\int_{F^\times} v'(x)v''(x) \vert x \vert^{-1} d^\times x.
\end{equation}

\begin{lemma}\label{2SClemma}
We have $0\neq \Phi\in \Hom_B ( V_2 \otimes V_3, \chi_1^{-1}\delta^{\frac{1}{2}})$.
\end{lemma}
\noindent{\it Proof: } Since $v_2$ and $v_3$ are given by the  characteristic function of $\mathcal{O}^\times$, 
$\Phi(v_2,v_3)=1\neq 0$. 
By (\ref{central}), $\Phi$ respects the central action. Since $\psi\overline{\psi}=1$, 
 $\Phi$ is also equivariant with respect to the action of $N$. Finally, for any $a\in F^\times$, 
 $$\Phi(\left(\begin{smallmatrix} a & 0 \\ 0 & 1 \end{smallmatrix}\right)\cdot v',\left(\begin{smallmatrix} a & 0 \\ 0 & 1 \end{smallmatrix}\right)\cdot v'')=\int_{F^\times} v'(ax)v''(ax) \vert x \vert^{-1} d^\times x
=\vert a \vert\Phi(v',v'')=$$
$$=(\chi_1^{-1}\delta^{\frac{1}{2}}) \left(\begin{smallmatrix} a & 0 \\ 0 & 1 \end{smallmatrix}\right) \Phi(v',v'').$$
 \hfill $\square$

It follows then from \cite[Proposition 4.5.5]{B}  that for any $v\otimes v' \otimes v''\in  V_1 \otimes V_2 \otimes V_3$ we have 
\begin{equation}\label{formula-for-l}
\ell(v\otimes v'\otimes v'')=\int_K v(k)\Phi(k\cdot v',k\cdot v'')dk.
 \end{equation}
 
\subsection{The case of unequal conductors}
In this subsection we assume that $n_2\neq n_3$, so
$n_2< n_3$. Since $V_1$ is minimal, it 
follows then from (\ref{central}) that $n_1< n_3$.

We first show that $\gamma^{n_3-n_1}\cdt v_1\otimes v_2\otimes v_3$ is a test vector. Since $\Phi(v_2,v_3)\neq 0$ by lemma \ref{2SClemma}, it follows that   $0\neq\ell(\bullet \otimes v_2\otimes v_3)\in 
\widetilde{V_1}^{I_{n_3},\omega_1^{-1}}$, hence there 
exists $0\leq i\leq n_3-n_1$ such that 
$\ell(\gamma^{i}\cdt v_1\otimes v_2\otimes v_3)\neq 0$.
Now, for every $0\leq i< n_3-n_1$, we have  
$$I_{{n_3}-1}  \subset  \gamma^{i}I_{n_1}\gamma^{-i}\cap I_{n_2}$$ 
hence    
$$\ell(\gamma^{i} \cdt v_1\otimes v_2\otimes \bullet)\in 
\widetilde{V_3}^{I_{n_3-1},\omega_3^{-1}}=\{0\}.$$
Therefore  $\ell(\gamma^{n_3-n_1}\cdt v_1\otimes v_2\otimes v_3)\neq 0$
as wanted. 

Next, we show that $ v_1\otimes \gamma^{n_3-n_2}\cdt v_2\otimes v_3$ is a test vector, assuming that 
$\ell(\bullet \otimes \gamma^{n_3-n_2}\cdt v_2\otimes v_3)\neq 0$.
As in the previous paragraph, there exists $0\leq i\leq n_3-n_1$ such that 
$\ell(\gamma^{i}\cdt v_1\otimes \gamma^{n_3-n_2}v_2\otimes v_3)\neq 0$. 
Moreover,  for every $0<i\leq  n_3-n_1$, we have  
$$\gamma I_{{n_3}-1}\gamma^{-1} \subset  \gamma^{i}I_{n_1}\gamma^{-i}\cap \gamma^{n_3-n_2}I_{n_2}  \gamma^{n_2-n_3}$$ 
hence    
$$\ell(\gamma^{i} \cdt v_1\otimes \gamma^{n_2-n_3} v_2\otimes \bullet)\in 
\widetilde{V_3}^{\gamma I_{n_3-1}\gamma^{-1},\omega_3^{-1}}=\{0\}.$$
 Therefore  $\ell( v_1\otimes \gamma^{n_3-n_2}\cdt v_2\otimes v_3)\neq 0$
as wanted. 

Finally, we prove the above assumption that $\ell(\bullet \otimes
\gamma^{n_3-n_2}\cdt v_2\otimes v_3)\neq 0$. 

Recall that $\left(\begin{smallmatrix} 0  & 1 
\\ \pi^{n_i} &  0\end{smallmatrix}\right)\cdot v_i$ is sent by the isomorphism $ V_i\otimes \omega_i^{-1}\cong\widetilde{V_i}$ to 
a new vector in $\widetilde{V_i}$. Moreover by (\ref{twist})
any test vector in $\widetilde{V_1}\otimes \widetilde{V_2}\otimes
\widetilde{V_3}$ yields a test vector in $V_1\otimes V_2 \otimes V_3$.  
By applying  lemma \ref{2SClemma}
to $\widetilde{V_1}\otimes \widetilde{V_2}\otimes
\widetilde{V_3}$ one gets 
$$\ell\left(\bullet \otimes \left(\begin{smallmatrix} 0  & 1
\\ \pi^{n_2} &  0\end{smallmatrix}\right)\cdot v_2\otimes \left(\begin{smallmatrix} 0  & 1
\\ \pi^{n_3} &  0\end{smallmatrix}\right)\cdot v_3\right)\neq 0,$$
hence $\ell(\bullet \otimes \gamma^{n_3-n_2}\cdt v_2\otimes v_3)\neq 0$. This completes the proof of theorem \ref{vt-main} in this case.

\subsection{The case of equal conductors}

In this subsection we assume that $n_2=n_3$, hence  $V_1$ is a 
ramified principal series. Since $V_1$ is minimal, it follows then from (\ref{central}) that $n_1< n_3$. By (\ref{formula-for-l}) and lemma 
\ref{calcul-SP1} we have 
$$\ell(\gamma^{n_3-n_1}\cdt v_1\otimes v_2\otimes v_3)=
\int_{I_{n_3}} v_1(k)\Phi(k\cdot v_2,k\cdot v_3)dk=
{\alpha_1}^{n_3-n_1}\int_{I_{n_3}} (\omega_1\omega_2\omega_3)(d) dk \neq 0.$$
where  $d$ is the lower right coefficient of  $k$.

Recall again that $\left(\begin{smallmatrix} 0  & 1 
\\ \pi^{n_i} &  0\end{smallmatrix}\right)\cdot v_i$ is sent by the isomorphism $ V_i\otimes \omega_i^{-1}\cong\widetilde{V_i}$ to 
a new vector in $\widetilde{V_i}$.  Moreover by (\ref{twist})
any test vector in $\widetilde{V_1}\otimes \widetilde{V_2}\otimes
\widetilde{V_3}$ yields a test vector in $V_1\otimes V_2 \otimes V_3$, 
hence 
$$\ell\left(\omega_1^{-1}(\det(\gamma^{n_3-n_1}))\gamma^{n_3-n_1}\cdot 
\left(\begin{smallmatrix} 0  & 1\\ \pi^{n_1} &  0\end{smallmatrix}\right)
\cdot v_1 \otimes \left(\begin{smallmatrix} 0  & 1
\\ \pi^{n_2} &  0\end{smallmatrix}\right)\cdot v_2\otimes \left(\begin{smallmatrix} 0  & 1
\\ \pi^{n_3} &  0\end{smallmatrix}\right)\cdot v_3 \right)\neq 0,$$
$$\ell(v_1\otimes v_2\otimes v_3)=\omega_1(\pi^{n_3-n_1})
\ell\left( \left(\begin{smallmatrix} 0  & 1
\\ \pi^{n_3} &  0\end{smallmatrix}\right)^{-1}
\gamma^{n_3-n_1}\left(\begin{smallmatrix} 0  & 1
\\ \pi^{n_1} &  0\end{smallmatrix}\right)\cdot v_1
\otimes v_2\otimes v_3\right)\neq 0.$$

This completes the proof of Theorem \ref{vt-main}. \hfill $\square$

\section{Test vectors in  reducible induced  representation} \label{reductible}

In this section, we generalize the local part of the paper \cite{HS} by Michael
Harris and Anthony Scholl on trilinear forms and test vectors when 
some of the $V_i$'s are  reducible principal series of $G$.  
The results of Harris and Scholl have as a global application the fact 
that a certain subspace,
constructed by Beilinson, in the motivic cohomology of the product of
two modular curves is a line. However, we are not going to follow them
in this direction.

As in \cite{HS}, we will only consider 
reducible principal series having infinite dimensional subspaces 
(see section  \ref{sous-espace}), since  for those having infinite
dimensional quotients  (see section \ref{sous-quotient}) test vector
can be obtained by preimage of  test vectors in the quotient. 
It follows then from \cite[Propositions 1.5, 1.6  and 1.7]{HS}  that under  
the assumption (\ref{central}): 
\begin{equation}
\dim\Hom_G ( V_1 \otimes V_2 \otimes V_3, \C)=1.
\end{equation}
This is particularly interesting for $V_1=V_2=V_3=
 \Ind_{B}^{G} (\delta^{\frac{1}{2}})$ since, according to
Theorem \ref{vt-000}(ii),  the space  $\Hom_G ( \mathrm{St} \otimes \mathrm{St}
\otimes \mathrm{St} , \C)$ vanishes. 

\begin{remark}
The case when for $1\leq i\leq 3$, $V_i=\Ind_{B}^{G}((\eta_i\circ\det)
  {\delta}^{\frac{1}{2}})$, with  $\eta_1 \eta_2 \eta_3$ non-trivial (quadratic),
   is not contained explicitly in \cite{HS}, but can 
be handled as follows. Since
$$\Hom_G \Bigl(\Ind_{B}^{G} ((\eta_1 \eta_2\circ\det) \delta^{\frac{1}{2}}), 
\Ind_{B}^{G} (({\eta_3}^{-1}\circ\det) {\delta}^{-\frac{1}{2}}) \Bigr) =0,$$
it follows easily from the short exact sequence
(\ref{speciale-sous-espace}) for $V_3$ that  there is  an isomorphism
$$\Hom_G ( V_1  \otimes V_2,  \widetilde{V_3} )\xrightarrow{\sim}
\Hom_G ( V_1  \otimes V_2,  \mathrm{St}\otimes\eta_3^{-1} ),$$
and the latter space is one dimensional by \cite[Proposition 1.6]{HS}.
\end{remark}

In \cite{HS}, Harris and Scholl also exhibit test vectors when the three
representations involved have a line of $K$-invariant vectors. 
The following proposition generalizes their results.

\begin{proposition} \label{vt-reductible} 
\begin{enumerate}
\item Suppose that for $1\leq i\leq 3$, $V_i=\Ind_{B}^{G}((\eta_i\circ\det)
  {\delta}^{\frac{1}{2}})$, with $\eta_i$  unramified character   such
  that $\eta_1^2 \eta_2^2 \eta_3^2=1$.  
 Then $v_1^K \otimes v_2^K \otimes v_3^K$ is a test vector.

\item  Suppose that for $1\leq i\leq 2$, $V_i=\Ind_{B}^{G}((\eta_i \circ\det)
  {\delta}^{\frac{1}{2}})$, with $\eta_i$  unramified, and 
$V_3$ is irreducible  such that 
$\eta_1^2 \eta_2^2\omega_3=1$.
Then $\gamma^{n_3} \cdt v_1^K \otimes v_2^K \otimes v_3$ and 
$v_1^K \otimes \gamma^{n_3} \cdt v_2^K \otimes v_3$ are test vectors.

\item  Suppose that $V_1=\Ind_{B}^{G}((\eta_1\circ\det)
   {\delta}^{\frac{1}{2}})$ with $\eta_1$  unramified, and that 
   $V_2$ and $V_3$ are irreducible  with $\eta_1^2\omega_2\omega_3=1$. 
Suppose that either $V_2$ is  non-supercuspidal and minimal, or
$V_2$ and $V_3$ are both supercuspidal with distinct 
conductors.  Then exactly one of the following holds: 
\begin{enumerate}
\item   $n_3 >n_2$ and  $v_1^K \otimes\gamma^{n_3-n_2} \cdt v_2\otimes  v_3$ 
and $\gamma^{n_3} \cdt v_1^K \otimes v_2\otimes  v_3$ are both test vectors;
\item $n_3=n_2$ and , for every $i$, $0 \leq i \leq n_3$, $
  \gamma^i\cdt v_1^K \otimes v_2 \otimes v_3$  is a test vector; 
\item  $V_2$ is special, $n_3=0$, and $v_1^K \otimes v_2\otimes  \gamma \cdt v_3$  
and $ \gamma \cdt v_1^K \otimes v_2\otimes v_3$ are both test vectors. 
\end{enumerate} 
\end{enumerate} 
\end{proposition}

\begin{remark}
One should observe that the test vectors in Proposition \ref{vt-reductible} :

$\relbar$ do not belong to {\it any  proper 
  subrepresentation} of $V_1\otimes V_2\otimes V_3$;

$\relbar$ are fixed by larger open compact subgroups of $G\times G
\times G$, than those fixing the  test vectors in the  
irreducible subrepresentation of $V_1\otimes V_2\otimes V_3$ given by 
 Theorem \ref{vt-main}. 
\end{remark}

\noindent{\it Proof: }  
As explained in the introduction,  twisting allows us to assume that 
$\eta_1=\eta_2=1$.

\noindent{\it  (i) } If $\eta_3=1$ this is \cite[Proposition 1.7]{HS}. 
Otherwise $\eta_3$ is the unramified quadratic character and 
we consider  Prasad's short exact sequence (\ref{courtesuite}): 
\begin{equation}
0 \rightarrow {\rm ind}_{T}^{G} 1  
  \xrightarrow{\mathrm{ext}} V_1\otimes V_2
  \xrightarrow{\mathrm{res}} \Ind_{B}^{G}  {\delta}^{\frac{3}{2}}\rightarrow 0.
\end{equation}
Since $\Hom_G (\Ind_{B}^{G} \delta^{\frac{3}{2}}, 
\Ind_{B}^{G} (({\eta_3}^{-1}\circ\det) {\delta}^{-\frac{1}{2}}) ) =0$, one
has : 
$$\Hom_G ( V_1  \otimes V_2,  \widetilde{V_3} )\xrightarrow{\sim}
\Hom_G ( {\rm ind}_{T}^{G}1 , \widetilde{V_3}  )
\xrightarrow{\sim} \Hom_T ( 1 ,\widetilde{V_{3\vert T}}). $$
Denote by  $\varphi$  a generator of the latter. It follows from  the proof
of Lemma \ref{lemmeV3},  where  the irreducibility of $V_3$ is not
used, only it's smoothness, that 
$\varphi(v_3^K)\neq 0$
(the point is that by (\ref{speciale-sous-espace}), a basis of the
$I_s$-invariants in $V_3$ is given by $\gamma^i \cdt v_3^K$ for $0\leq
i\leq s$).  

It follows then by exactly the same argument as in the proof of 
 \cite[Theorem 5.10]{P}, that $v_1^K \otimes v_2^K
\otimes v_3^K$ is a test vector. The only point  to check is that the 
denominator in the formula displayed in the middle of \cite[page
20]{P} does not vanish. 

\bigskip
\noindent{\it  (ii) }
For $n_3 =0$, this is \cite[Proposition 1.6]{HS}.

For  $n_3 \geq 1$, again by  Lemma \ref{lemmeV3} we have
$\varphi(v_3^K)\neq 0$ and the usual process, as in the proof of 
Theorem \ref{vt-00n}, allows to prove that $\gamma^{n_3}\cdt v_1^K \otimes
v_2^K \otimes v_3$ and  
$ v_1^K \otimes \gamma^{n_3}\cdt v_2^K \otimes v_3$ are test vector. 

\bigskip
\noindent{\it (iii)(a) }
If $V_2$ and $V_3$ are both supercuspidal the claim follows
from lemma \ref{inclusion} by exactly  same arguments that allowed to 
prove theorem \ref{vt-main} in this case. So we can assume that 
$V_2$ is non-supercuspidal and minimal. 

First we choose a model of $V_2$ such that   $\mu_2$ is unramified and
consider the exact sequence (\ref{courtesuite}): 
$$ 0 \rightarrow {\rm ind}_{T}^{G}(
{\delta}^{-\frac{1}{2}} \chi_2)   
  \xrightarrow{\mathrm{ext}} V_1 \otimes V_2 
  \xrightarrow{\mathrm{res}} \Ind_{B}^{G} ( {\delta}\chi_2)\rightarrow 0. $$

If $\Hom_G \Bigl( \Ind_{B}^{G} ({\delta}\chi_2 ),\widetilde{V_3} \Bigr)=0$,
then we obtain  isomorphisms
$$\Hom_G \Bigl( V_1 \otimes V_2 , \widetilde{V_3} \Bigr)
\xrightarrow{\sim}
\Hom_G \Bigl({\rm ind}_{T}^{G}( 
{\delta}^{-\frac{1}{2}} \chi_2)  , \widetilde{V_3} \Bigr)
\xrightarrow{\sim}
\Hom_T \Bigl({\delta}^{\frac{1}{2}}\chi_2
, \widetilde{V_{3\vert T}} \Bigr)$$ 
and as in  section \ref{Preuve-main} we obtain that
$ v_1^K \otimes \gamma^{n_3-n_2}\cdt v_2 \otimes v_3$ is a test vector. 

\medskip

If $\Hom_G \Bigl( \Ind_{B}^{G} ({\delta}\chi_2 
),\widetilde{V_3} \Bigr)\neq 0$, then $n_3>n_2$ implies that 
there exists an unramified character $\eta$ such that 
${\delta}\chi_2 =(\eta\circ\det) \delta^{-\frac{1}{2}}$ 
 and $\widetilde{V_3}=\eta \otimes \mathrm{St}$. So $n_2=0$ and $n_3=1$. 
It is easy then to check that the image of 
$ v_1^K\otimes \gamma \cdt v_2 \otimes \in V_1 \otimes V_2$ by $\mathrm{res}$ 
is not a  multiple of $\eta \circ \det$, 
hence it yields a non zero element of $\widetilde{V_3}$. 
Since $\gamma^{-1}K\gamma \cap K =I$, it is actually a non zero element of 
$\widetilde{V_3}^{I,\omega_3^{-1}}$, hence 
$v_1^K\otimes \gamma \cdt v_2 \otimes  v_3$ is a test vector.

\medskip
By choosing a model of $V_2$ with $\mu'_2$ 
unramified, and 
applying the above arguments  to 
 $V_2  \otimes \Ind_{B}^{G} ({\delta}^{\frac{1}{2}})$ 
one can prove that 
$\gamma^{n_3} \cdt v_1^K \otimes  v_2 \otimes v_3$ is a test vector.

\bigskip
\noindent{\it (iii)(b)} For $n_3 =0$, this is \cite[Proposition 1.5]{HS}.

For  $n_3 \geq 1$, assume first that $\Hom_G \Bigl(V_2,
\widetilde{V_3}\Bigr)\neq 0$. Then the $G$-invariant trilinear form 
on $V_1\otimes V_2 \otimes V_3$ is obtained by composing 
$\mathrm{proj}^*_1\otimes \mathrm{id}\otimes\mathrm{id} $ with the natural pairing between $V_2\simeq \widetilde{V_3}$ and $V_3$. Since the natural
pairing between $\widetilde{V_3}$ and $V_3$ is non-zero on a couple of
new vectors, it follows that for all $i$, 
$\gamma^i \cdt v_1^K  \otimes v_2 \otimes  v_3$ is a test vector.

If  $\Hom_G \Bigl(V_2,  
\widetilde{V_3}\Bigr)=0$, we apply  the  techniques of section
\ref{Preuve-main} to $V_2 \otimes V_3 \otimes \Ind_{B}^{G} (
{\delta}^{\frac{1}{2}})$. There are  isomorphisms
$${\rm Hom}_G  ( V_2 \otimes V_3 \otimes V_1 , \C ) \xrightarrow{\sim} 
{\rm Hom}_T \Bigl( \chi_3\chi'_2  , {\Ind_{B}^{G}
  ({\delta}^{-\frac{1}{2}})_{\vert T}} \Bigr) \xrightarrow{\sim} 
{\rm Hom}_T \Bigl( \chi_3\chi'_2  , \widetilde{V_{1\vert T}} \Bigr).$$
Taking a generator $\varphi$ of the latter space, one has
$\varphi(\gamma^i \cdt v_1^K) \not=0$ for all $i$, by adapting the proof 
of  Lemma \ref{lemmeV3} as above. Then exactly the same computations
as in the proof of Theorem \ref{vt-main}(b) show that  $
\gamma^i v_1^K \otimes v_2 \otimes v_3$ is a test vector, for all $0\leq i\leq
n_3$.

\noindent{\it (iii)(c) } This case follows from (iii)(a) applied to 
$ V_1 \otimes V_3 \otimes V_2$.  \hfill $\Box$

\end{document}